\documentclass[11pt,leqno,twoside]{amsart}
\usepackage{amsmath,amsthm,amsfonts}
\usepackage{amssymb,verbatim,url,graphicx}
\usepackage[percent]{overpic}

\graphicspath{{figs/}}
\newcommand{\figr}[1]{Figure~\ref{fig:#1}}
\newcommand{\fig}[5]{%
\begin{figure}
\begin{overpic}[width=#2]{#1}
#4
\end{overpic}
\caption[#3]{#5}
\label{fig:#1}
\end{figure}%
}
\newcommand{\figtwo}[7]{%
\begin{figure}
\begin{overpic}[width=#3]{#1} #5 \end{overpic}
\hspace{1cm}
\begin{overpic}[width=#3]{#2} #6 \end{overpic}
\caption[#4]{#7} \label{fig:#1-#2}
\end{figure}%
}

\newtheorem{theorem}{Theorem}[section]
\newtheorem*{maintheorem}{Main Theorem}
\newtheorem{lemma}[theorem]{Lemma}
\newtheorem{proposition}[theorem]{Proposition}
\newtheorem{corollary}[theorem]{Corollary}
\newtheorem*{maincorollary}{Corollary}

\newcommand{\prop}[1]{Proposition~\ref{pr:#1}}
\newcommand{\cor}[1]{Corollary~\ref{co:#1}}

\theoremstyle{definition}
  \newtheorem{definition}[theorem]{Definition}
  
\theoremstyle{remark}
  \newtheorem*{remark}{Remark}

  \newtheorem*{example}{Example}

\newcommand{\ia}{{\rm({\it i\/}\rm)}}
\newcommand{\ii}{{\rm({\it ii\/}\rm)}}

\newcommand{\area}{\operatorname{area}}
\def\alex-emb/{Alexandrov-embedded}
\newcommand{\after}{\circ}
\newcommand{\C}{\mathbb{C}}
\newcommand{\cmc}{{\footnotesize{CMC}}~}
\def\CMC/{{\footnotesize{CMC}}}

\def\d{\partial}
\newcommand{\dist}{\operatorname{dist}}
\newcommand{\D}{\mathcal{D}}
\newcommand{\Dk}{\D_k}
\newcommand{\Dbox}{{D^{{}_\Box}_{\vphantom{j}}}}  

\newcommand{\E}{\mathcal{E}}

\renewcommand{\epsilon}{\varepsilon}
\newcommand{\grad}{\mathrm{grad}}
\renewcommand{\H}{\mathbb{H}}

\newcommand{\ident}{\equiv}
\newcommand{\isom}{\cong}
\renewcommand{\Im}{\operatorname{Im}}

\let\dotlessi\i
\renewcommand{\i}{\mathbf{i}}  
\renewcommand{\k}{\mathbf{k}}  
\newcommand{\K}{ K} 

\newcommand{\M}{\mathcal{M}}
\newcommand{\Mk}{\mathcal{M}_k}
\newcommand{\N}{\mathbb{N}}
\let\paragraph=\S  

\renewcommand{\P}{{\mathcal P}}
\newcommand{\Pk}{\P_k}

\renewcommand{\phi}{\varphi}
\newcommand{\Pik}{\Pi_{\k}}               
\def\wtilde#1{{\mskip4mu\widetilde{\mskip-4mu #1\mskip-1mu}\mskip1mu}}
\def\mtilde#1{{\mskip7mu\widetilde{\mskip-7mu #1\mskip-7mu}\mskip7mu}}
\def\plustilde#1{\wtilde{#1}{}^+}
\newcommand{\R}{\mathbb{R}}

\renewcommand{\S}{\mathbb{S}}
\newcommand{\sfrac}[2]{#1/#2} 
\newcommand{\SO}{\mathsf{SO}}
\newcommand{\setm}{\smallsetminus}
\newcommand{\T}{\mathcal{T}}

\newcommand{\Tk}{{\mathcal{T}}_k}
\renewcommand{\theta}{\vartheta}
\def\takes{\colon}
\def\uhalf/{upper half}

\newcommand{\X}{\mathcal{X}}
\newcommand{\Y}{\mathcal{Y}}
\newcommand{\Z}{\mathbb{Z}}
\newcommand{\Zodd}{\mathbb{Z}_\mathrm{odd}} 

\let\mgp=\marginpar \marginparwidth18mm \marginparsep1mm
\def\marginpar#1{\mgp{\raggedright\tiny #1}}

\begin{document}

\title{Coplanar constant mean curvature surfaces}
\date{2007 December 5}

\author[Grosse-Brauckmann]{Karsten~Gro{\ss}e-Brauckmann}
\address{Technische Universit\"at Darmstadt, Fb.~Mathematik (AG~3),
    Schlossgartenstr.~7, 64289~Darmstadt, Germany}
\email{kgb@mathematik.tu-darmstadt.de}

\author[Kusner]{\\Robert~B.~Kusner}
\address{Mathematics Dept., University of Massachusetts,
      Amherst~MA~01003, USA}
\email{kusner@math.umass.edu}

\author[Sullivan]{John~M.~Sullivan}
\address{Technische Universit\"at Berlin, MA 3--2,
Str.~des 17.~Juni~136, 10623~Berlin, Germany}
\email{sullivan@math.tu-berlin.de}

\begin{abstract}
  We consider constant mean curvature surfaces with finite topology,
  properly embedded in three-space in the sense of Alexandrov.
  Such surfaces with three ends and genus zero were constructed
  and completely classified by the authors~\cite{gkspnas, gks}.  
  Here we extend the arguments to the case of an arbitrary number of ends,
  under the assumption that the asymptotic axes of the
  ends lie in a common plane:
  we construct and classify the entire family of these genus-zero,
  coplanar constant mean curvature surfaces.
\end{abstract}
\maketitle

\begin{center}
\small \emph{Dedicated to Hermann Karcher\\
on the occasion of his sixty-fifth birthday.}
\end{center}
\bigskip

\section*{Introduction}

The search for examples of surfaces with constant mean curvature
dates back to the nineteenth century.
Unlike minimal surfaces, which have a Weierstrass representation,
the case of nonzero constant mean curvature (\CMC/) does not readily yield 
explicit examples; 
it thus appeared that constant mean curvature surfaces were scarce.  
In recent decades, however, several construction methods have evolved.
Wente's discovery of an immersed torus of constant mean curvature
initiated the integrable-systems approach explored 
by Abresch, Bobenko, Pinkall, Sterling and others; 
this in turn led to a generalized Weierstrass representation 
by Dorfmeister, Pedit and Wu.
The ground-breaking analytic perturbation (or gluing) techniques of Kapouleas 
have been extended by Mazzeo, Pacard, Pollack, Ratzkin and others.  
Conjugate surface techniques, introduced by Lawson, 
were developed further by Karcher 
and---combined with new analytic methods---by the present authors.  
An abundance of known \CMC/ surfaces has emerged, 
inviting their systematic study and classification.

With the conjugate surface method, we can study \CMC/
surfaces with sufficient symmetries.  This method is powerful in
that it enables us to study entire families of surfaces, that is,
components of suitably defined moduli spaces.
In contrast, the gluing techniques mentioned above, while more widely
applicable, tend to produce examples only near a boundary of the
relevant moduli space.
The generalized Weierstrass representation also has the potential to deal
with moduli spaces, even without any symmetry assumptions,
but it has difficulty detecting key geometric features
such as embeddedness.

The \emph{Delaunay unduloids}~\cite{del} are
the noncompact embedded \CMC/ surfaces
of revolution.  These singly periodic surfaces decompose along
their \emph{necks} (minimizing closed geodesics) into \emph{bubbles}.  
All known embedded \CMC/ surfaces display similar bubbles.
In light of experience with
doubly and triply periodic surfaces \cite{kar,law,kgb},
the discovery that bubbles can be created and deleted 
continuously on surfaces of finite topology~\cite{gk} was surprising.
This insight is what led us to expect our moduli spaces to be connected,
and thus to look for the existence proof presented here, 
based on the continuity method of~\cite{gks}.

Let us now make the setting of our paper precise.
A surface has \emph{finite topology} when it is homeomorphic to
a compact surface~$\Sigma$ of genus~$g$ from which a finite number~$k$ of 
points has been removed;  
a neighborhood of any of these punctures is called an \emph{end}.
If the surface has nonzero constant mean curvature, 
then by a scaling we may assume the mean curvature is identically~$1$; 
from now on \cmc will denote this condition.
As in \cite{gks} we consider a class of surfaces more general
than the embedded ones, namely those bounding an immersed solid.
To be precise, we define a finite topology \CMC/ surface~$M$
to be \emph{\alex-emb/} if~$M$ is properly immersed, if each end of~$M$ is
embedded, and if there exists a compact 3-manifold~$W$ with boundary
$\d W=\Sigma$, and a proper immersion
$F\takes W\setm\{q_1,\ldots,q_k\}\to\R^3$ whose boundary restriction
$f\takes \Sigma\setm\{q_1,\ldots,q_k\} \to\R^3$ parametrizes~$M$.  
Moreover, we assume that the (inward) mean-curvature normal~$\nu$ of~$M$ 
points into~$W$.  

A Delaunay unduloid is an embedded \CMC/ surface with two ends and genus zero;
it is determined by its \emph{necksize} $n\in(0,\pi]$, the circumference of any neck.
More generally, 
we define a \emph{$k$-unduloid} to be an \alex-emb/ \CMC/ surface with
finite genus and a finite number $k\in\N$ of ends.  Our terminology
is justified by the asymptotics result of Korevaar, Kusner and Solomon:
\begin{theorem}[\cite{kks}]  \label{th:kksasym}
  Each end of a $k$-unduloid is exponentially asymptotic to 
  a Delaunay unduloid, and thus has an asymptotic axis and 
  necksize $n\in(0,\pi]$.
\end{theorem}
For $k\le2$, all $k$-unduloids have genus zero, and their classification
is well-known: Alexandrov~\cite{al1,al2} showed the round sphere
is the only $0$-unduloid; Meeks~\cite{mee} showed there are
no $1$-unduloids; and the Delaunay unduloids themselves
are the only $2$-undu\-loids~\cite{kks}. 

A \cmc surface is called \emph{coplanar} 
if it is contained in a half-space of~$\R^3$. 
For $k\le3$, all $k$-unduloids are coplanar~\cite{mee}.
The Alexandrov reflection method extends to coplanar $k$-unduloids~\cite{kks},
showing that any such surface~$M$ is \emph{Alexandrov-symmetric} 
in the following sense:
it has a reflection plane~$P$    
such that $M\setm P$ has two connected
components $M^{\pm}$ whose Gauss images~$\nu(M^{\pm})$ are
contained in disjoint open hemispheres; 
moreover $\nu(M\cap P)$ is contained in the equator.
%
\begin{theorem}[\cite{kks}]  \label{th:kkscopl}
  A coplanar $k$-unduloid $M$ is contained in a slab 
  and is Alexandrov-symmetric across a plane~$P$ within that slab.  
  The plane~$P$ contains the asymptotic axes of the ends of~$M$.
\end{theorem}

We want to study moduli spaces of coplanar $k$-unduloids with genus zero.
Thus we identify two such $k$-unduloids 
when they differ by a rigid motion in space and we 
label the ends to avoid orbifold points at surfaces with extra symmetry.
It is convenient to represent any coplanar $k$-unduloid
with the horizontal $xy$-plane as its fixed mirror plane~$P$,
such that the labels of the ends occur in cyclic order:
\begin{definition}
  Let us fix $k\ge0$ points $q_1,\ldots,q_k$
  in cyclic order along the equator of~$\S^2$.
  We consider proper \alex-emb/ \CMC/ immersions of 
  $\S^2\setm \{q_1,\ldots ,q_k\}$ into~$\R^3$
  which are equivariant under reflection across the $xy$-plane.
  The \emph{moduli space}~$\Mk$ of {coplanar $k$-unduloids} 
  of genus zero consists of all such immersions, modulo 
  diffeomorphisms of the domain fixing the points~$q_i$,
  and modulo isometries of~$\R^3$.
  Uniform convergence in Hausdorff distance
  on compact subsets of~$\R^3$ defines the topology on~$\Mk$.
\end{definition}


Our main theorem classifies coplanar $k$-unduloids of genus zero in~$\Mk$.
This generalizes our classification of 
\emph{triunduloids} (with $k=3$) 
by triples of points in~$\S^2$ modulo rotation~\cite{gks},
which shows $\M_3$ is homeomorphic to the open three-ball.
For $k\ge4$, there is still a natural map from
coplanar $k$-unduloids of genus zero to $k$-tuples on~$\S^2$,
but this no longer suffices to classify the surfaces.
Instead, we need to employ a space~$\Dk$ of spherical metrics,
whose precise definition will be given below.
Our classification says that $\Mk$ is homeomorphic to~$\Dk$:
\begin{maintheorem}
  For each $k\ge3$, there is a homeomorphism $\Phi$
  from the moduli space~$\Mk$ of coplanar $k$-unduloids
  of genus~$0$
  to the connected \hbox{$(2k-3)$}--manifold~$\Dk$ of $k$-point metrics.  
  The necksizes of a $k$-unduloid $M\in\Mk$ are 
  the spherical distances of the $k$~consecutive pairs of
  completion points in~$\Phi(M)$.
\end{maintheorem}
As a consequence we have the following constraints 
on the necksizes of a coplanar $k$-unduloid of genus~$0$.
(These do not hold for higher genus, as shown by
the coplanar $k$-unduloids of genus~$1$ 
with all ends cylindrical~\cite{cylends}.)
\begin{maincorollary}
  Let $k\ge 3$.
  The necksizes $n_1,\ldots,n_k\in(0,\pi]$
  of any coplanar $k$-unduloid of genus~$0$ satisfy
  \begin{align}
    \label{necksizesumodd}
    && n_1+ \cdots+ n_k &\le (k-1)\pi &\text{for }& k \text{ odd}, &&\\
    \label{necksizesumeven}
    && n_1+\cdots + n_k & < k\pi      &\text{for }& k\ge4 \text{ even}, &&
  \end{align}
  and these inequalities are sharp.
  At most $k-2$ ends of the surface can be cylindrical,
  and examples of coplanar $k$-unduloids with $k-2$ cylindrical ends
  exist for each $k\ge 3$.
\end{maincorollary}

\begin{figure}
\hspace*{-7mm}
\begin{overpic}[width=.23\textwidth]{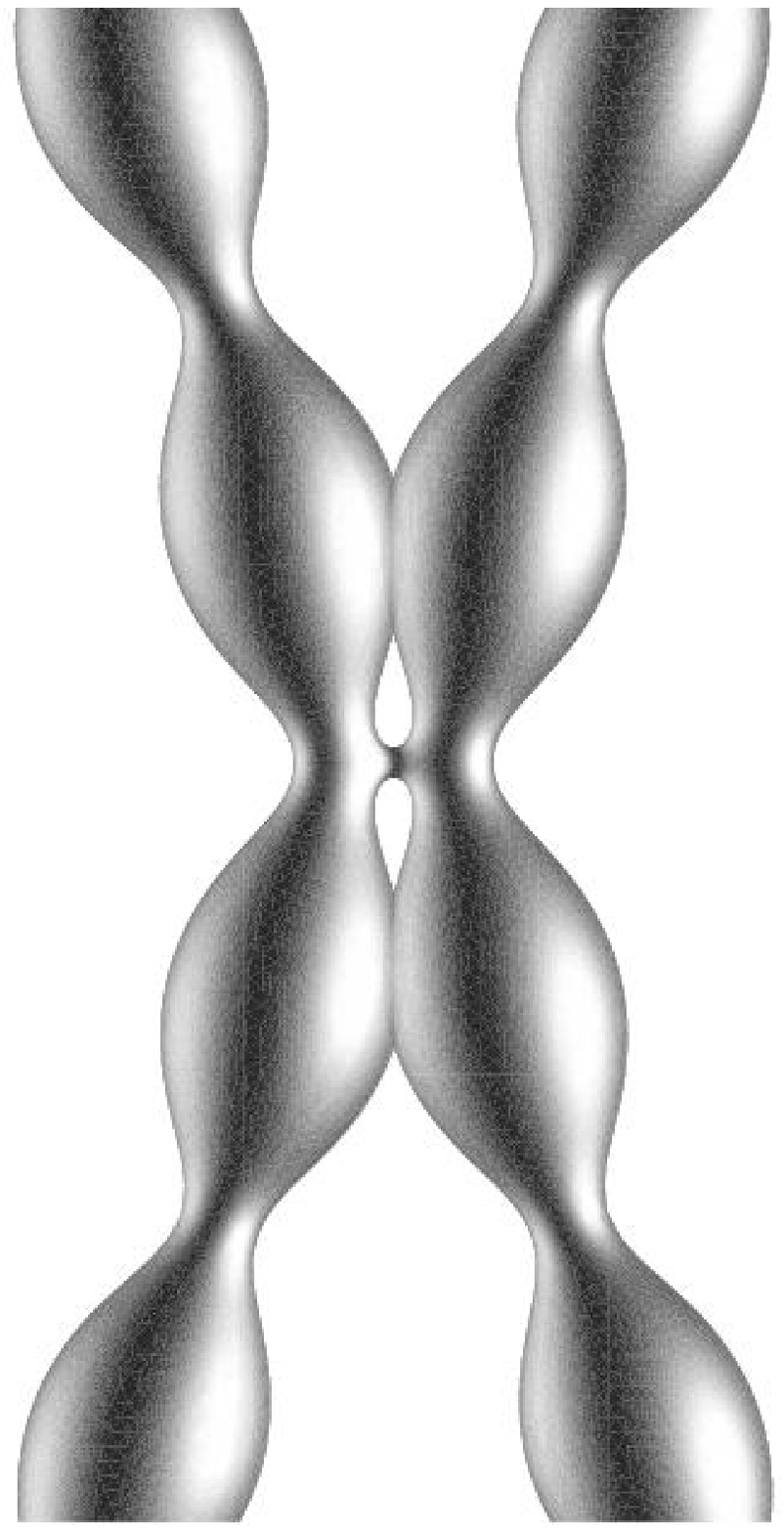} \end{overpic}
\hspace{-2mm}
\begin{overpic}[width=.23\textwidth]{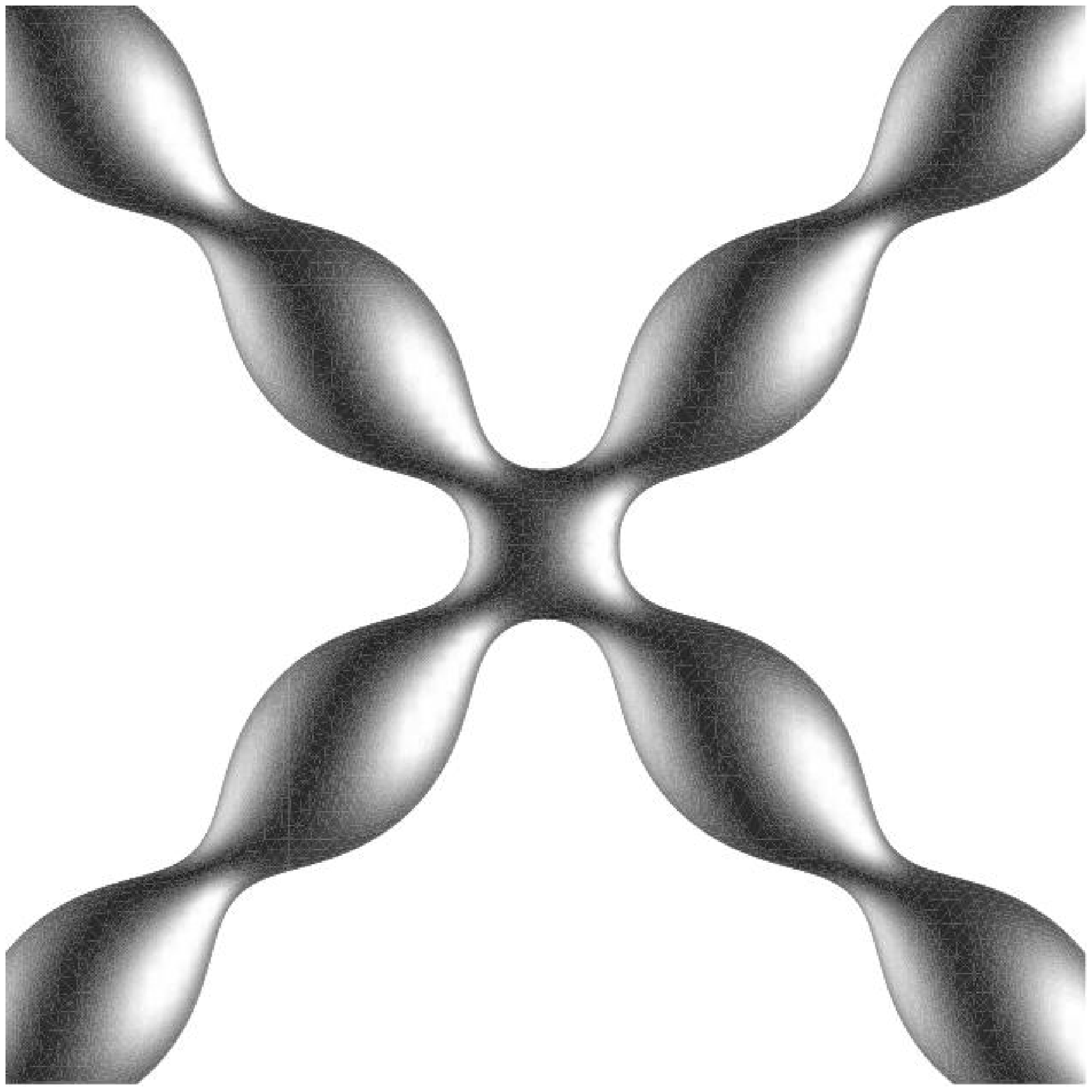} \end{overpic}
\hspace{5mm}
\begin{overpic}[width=.23\textwidth]{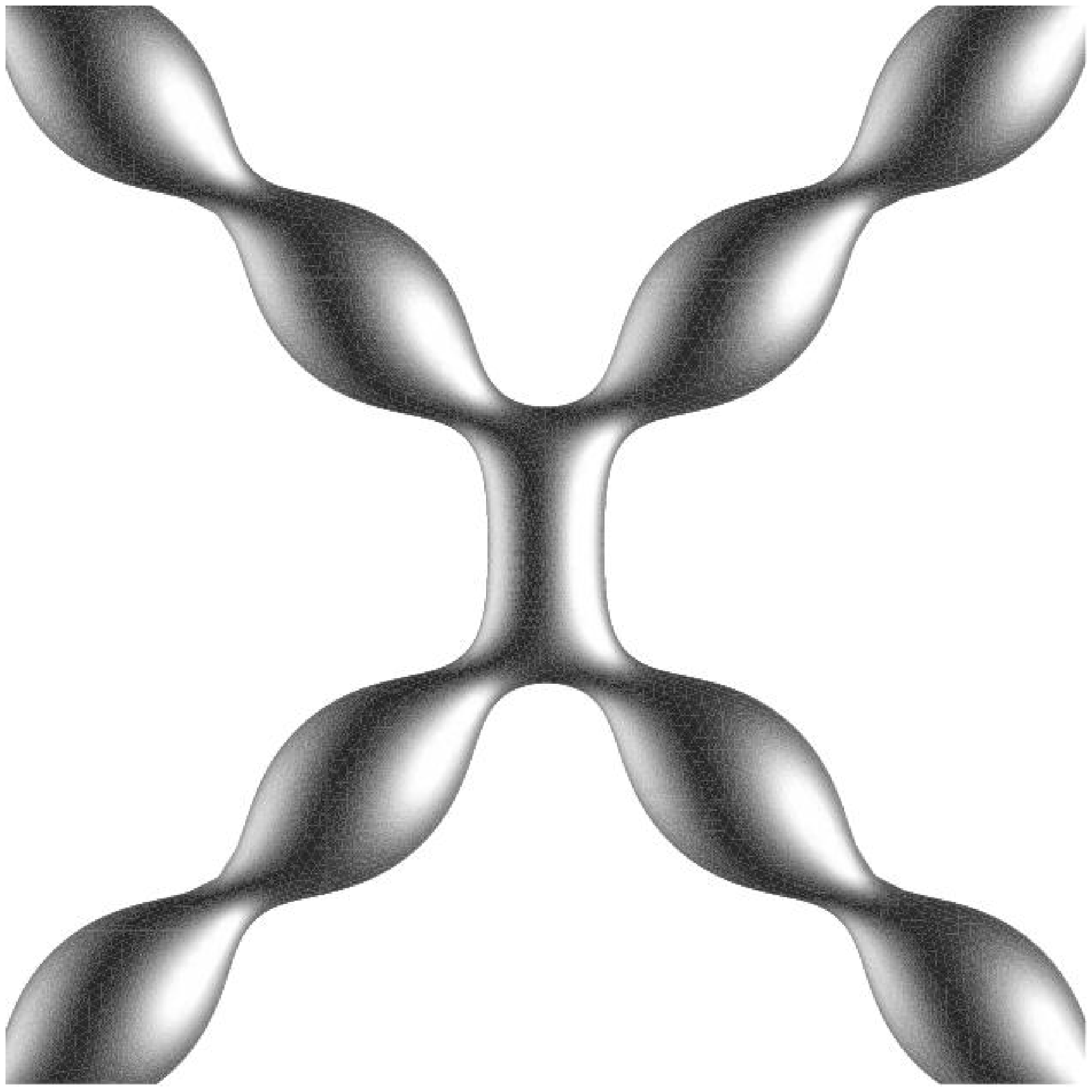} \end{overpic}
\hspace{2mm}
\begin{overpic}[width=.23\textwidth]{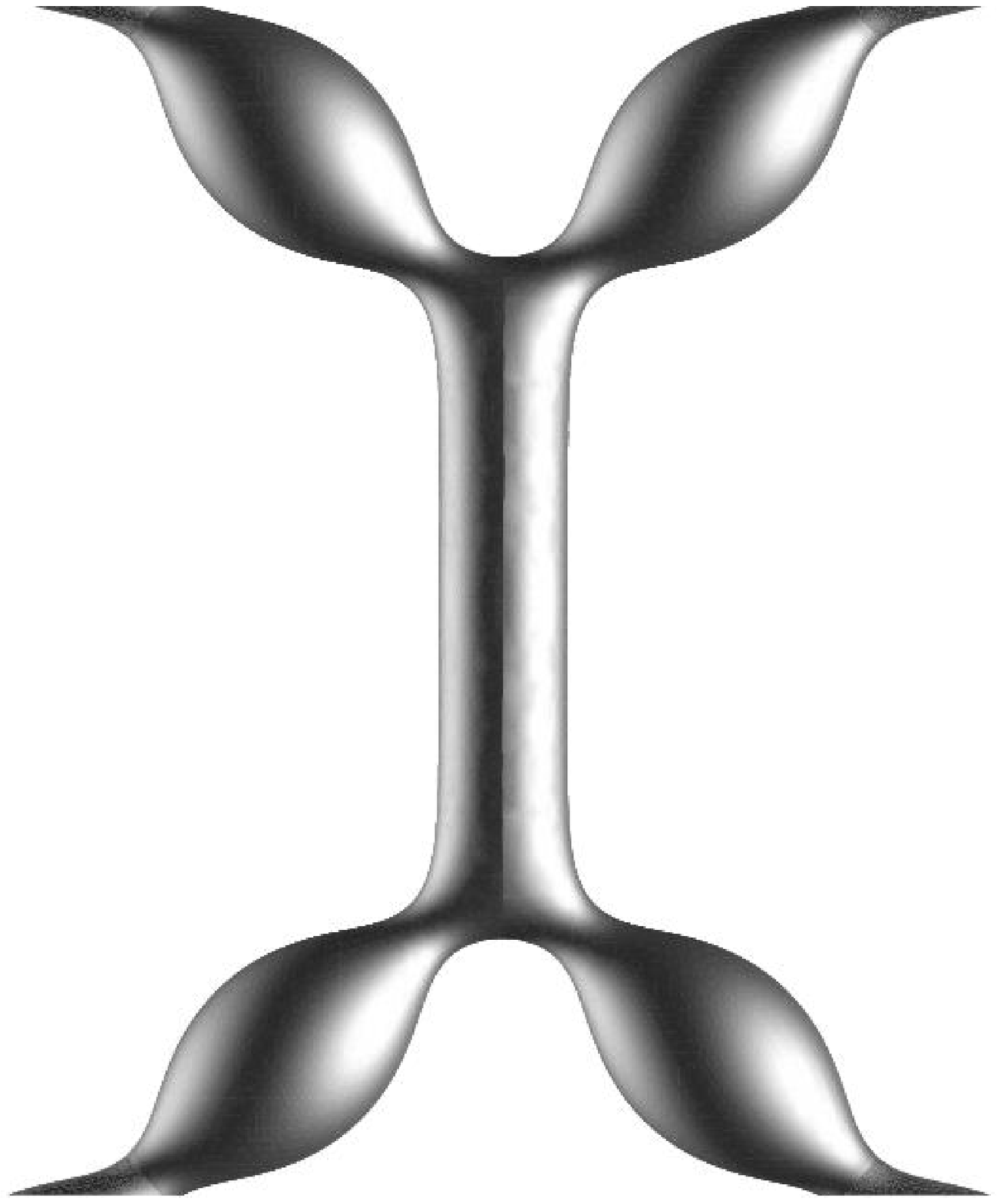} \end{overpic}
\hspace*{-5mm}
\caption[rhombic with cyl segment]{
  Our main theorem characterizes which sets of~$k$ asymptotic
  necksizes a genus-zero coplanar $k$-unduloid can attain.
  In addition, for $k=3$ the necksizes uniquely 
  determine the angles between the ends, as follows from force balancing.
  For $k\ge4$, however, there is no such dependence:
  the $4$-unduloids shown above all have all necksizes equal to~$\pi/2$,
  but the angles between the ends vary.  
}\label{fig:rhombic with cyl segment}
\end{figure}
 
Our theorem should be compared with the classification of
coplanar minimal surfaces by Cos{\'\dotlessi}n and Ros~\cite{cr}.
The classifying space for their minimal $k$-noids consists of planar polygons,
whose combinatorial description is well-known.
Another simplification arises because they are able
to prove directly that all their surfaces are nondegenerate.
We are not aware of any other classification results
for moduli spaces of minimal or \CMC/ surfaces with more than three ends.

Very recently, in joint work~\cite{GKKRS} with Korevaar and Ratzkin
(following up on~\cite{kkr}),
we have augmented the results of this paper as follows:
All coplanar $k$-unduloids of genus zero are nondegenerate,
so that $\Mk$ is a real-analytic manifold, and a tubular neighborhood
of it in the space of all $k$-unduloids is a $(3k-6)$--manifold.
Furthermore, our classifying map~$\Phi$ is a (real-analytic) diffeomorphism.
\medskip

We conclude this introduction with an overview of our paper.
In Section~1, we define the classifying map $\Phi$ as a composition
of the following maps:  By coplanarity and Theorem~\ref{th:kkscopl}, we first 
associate to each $M\in\Mk$ its \uhalf/~$M^+$. 
Since $M$ has genus zero, $M^+$ is an open disk.
It thus has a conjugate cousin, an isometric minimal
immersion to the three-sphere.
Hopf projection of this cousin disk to the $2$-sphere
produces a metric on $M^+$ locally isometric to $\S^2$, 
which we call the spherical $k$-point metric $\Phi(M)$.

Central for our paper is the analysis (in Sections~2 and~3)
of the space~$\Dk$ of spherical $k$-point metrics.
In Section~2,
we start with a combinatorial decomposition of each metric $D\in\Dk$.
In view of the classifying map, analyzing 
metrics amounts to analyzing coplanar $k$-unduloids,
yielding a decomposition of the surface into bubbles.
On cylindrical parts of a $k$-unduloid, such a decomposition 
is necessarily ambiguous.  For this reason, 
our decomposition cannot be unique, making the analysis somewhat involved.

In Section~3, we apply our combinatorial decomposition to
show that $\Dk$ is a connected manifold.
We analyze the example~$\D_4$ completely to find that 
that~$\D_4\isom \M_4$ is a five-ball, generalizing the
result of~\cite{gks} that $\D_3\isom\M_3$ is a three-ball.
Our combinatorial description is not explicit enough 
to decide whether~$\Dk$ is always a $(2k-3)$--ball.
(In ongoing joint work with Korevaar, using a M{\"o}bius-invariant
decomposition of $k$-point metrics, we expect to be
able to show that indeed $\Dk\isom\R^{2k-3}$ for all $k\ge3$.)

In Section~4, we prove injectivity and properness of our classifying map.
Injectivity follows by applying
a maximum principle to the cousin surfaces
(just as in the case $k=3$ of~\cite{gks}).
Properness involves \emph{a priori} estimates, 
corresponding to the closedness part of a traditional
continuity method.  The proof here is somewhat different than
in the case $k=3$ because of the possibility of small necks occurring
away from the ends.  Still, our characterization of compact subsets
of~$\Dk$ allows us to prove the needed area and curvature estimates.

In Section~5, we establish the surjectivity of our classifying map.
As in \cite{gks}, the crucial point is to apply a structure theorem
for real-analytic varieties.  Thus we must show
that our spaces~$\Mk$, consisting of coplanar surfaces
with a mirror symmetry, are real-analytic.
For this, we adapt the methods of~\cite{kmp},
which dealt with the space of all $k$-unduloids.
More generally, we show that for any fixed finite symmetry group,
the moduli space of symmetric $k$-unduloids is real-analytic,
a result which has applications beyond this paper~\cite{groisman}.
Then we can apply our version of the continuity method, as introduced
in~\cite{gks}, to prove the main theorem.

{\small
We gratefully acknowledge partial support from 
the Deutsche Forschungsgemeinschaft and the National Science Foundation.
Our figures of \cmc surfaces were produced with the 
Graphics Programming Environment \emph{Grape}; our
sketches of spherical polygons, with \emph{Spherical Easel}.
}

\section{The classifying map}       \label{se:class}

In this section we define the space~$\Dk$ of $k$-point spherical metrics
and the classifying map $\Phi\colon\Mk\to\Dk$.
We assume $k\ge 2$ to avoid trivial cases.

\subsection{Conjugate cousins and necksizes}

Consider a coplanar $k$-unduloid $M\in\Mk$;
it is Alexandrov symmetric across the $xy$-plane~$P$.
Take a component $M^+$ of $M\setminus P$
and orient $M^+$ such that the boundary labels occur in
increasing cyclic order:  
this means that $M^+$ lies to the left of~$\d M^+$.  
By a rotation of $M$ in $\R^3$, we may assume that
$M^+ := M\cap\{z>0\}$; we call $M^+$ the \emph{\uhalf/}.

When $M$ is embedded the ends occur in counterclockwise order in~$P$.
In other words, if at each end of~$M$ 
we truncate $P\cap M$ and join in a locally embedded way,
then we obtain an embedded loop in~$P$ with turning number one.

\begin{figure}
\begin{overpic}[width=.47\textwidth]{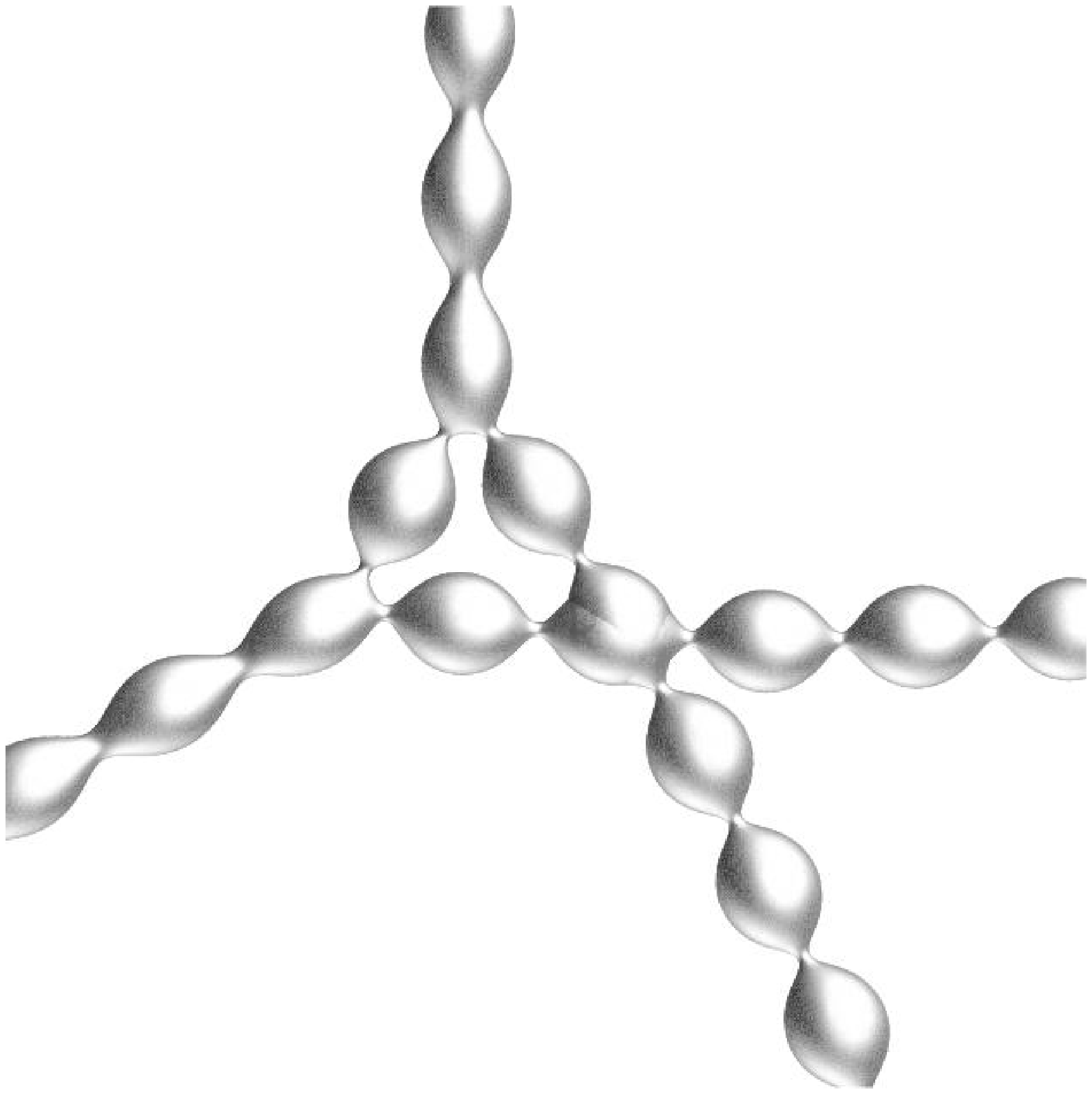} \end{overpic}
\hspace*{\fill}
\begin{overpic}[width=.47\textwidth]{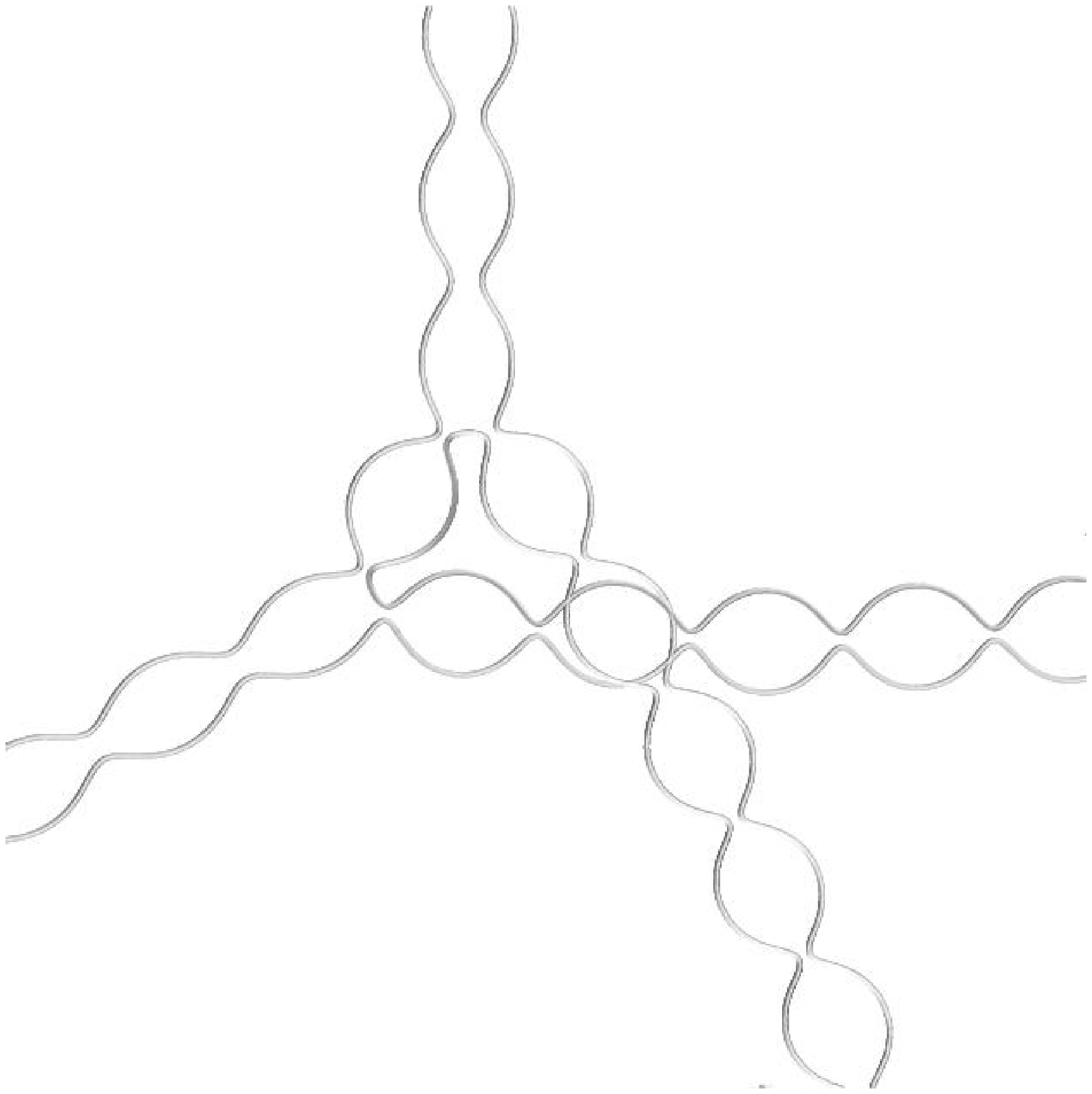} \end{overpic}
\caption{
  A coplanar 4-unduloid and its
  intersection with the symmetry plane~$P$.
  Typical coplanar $k$-unduloids with $k>3$ 
  have such crossed ends and are not embedded.
} 
\label{fig:crossed ends}
\end{figure}

However, in general ends of~$M$ may cross 
(as in~\figr{crossed ends}, right)
and thus no longer appear in the same order in~$P$.
We note that for large~$k$, most coplanar $k$-unduloids
have such crossed ends and thus fail to be embedded.
By Alexandrov symmetry, 
the projection of~$M^+$ to~$P$ gives an immersed disk,
so that the turning number of the corresponding boundary
loop is still one.

We now review the conjugate cousin relationship~\cite{law},
in the first-order form presented in \cite[Thm.~1.1]{gks}.
Given any simply connected \CMC/ surface, such as~$M^+$, there
is a cousin minimal surface $\plustilde M$ in the three-sphere~$\S^3$.
We regard~$\S^3$ as the Lie group of unit quaternions
and identify~$\R^3$ with its Lie algebra $T_1\S^3=\Im\H$, as in~\cite{gks}.
The cousin is determined uniquely up to left translation in~$\S^3$. 

The upper half~$M^+\subset\R^3$ is a disk isometric to its
minimal cousin $\plustilde M\subset\S^3$.
This isometry has the following effect on normals and tangent vectors.

If $\nu\in\S^2$ is the normal at $p\in M^+$,
then the normal $\tilde\nu$
at the corresponding point $\tilde p\in\plustilde M$
is simply the left translation of~$\nu$, that is, $\tilde \nu=\tilde p\nu$.
Left translation preserves the inner product, and so
for the unit quaternion $\k\in\S^2\subset\Im\H$ we find
$$
  \langle \tilde\nu, \tilde p\k\rangle_{\S^3} = \langle\nu, \k\rangle_{\R^3}.
$$
Since $\k$ is vertical,
Theorem~\ref{th:kkscopl} gives $\langle\nu,\k\rangle<0$ 
on the \uhalf/ $M^+$,  
and thus the cousin $\plustilde M$ is transverse to 
the \emph{$\k$-Hopf circles} $t\mapsto\tilde p(\cos t+\k\sin t)$. 

If $X$ is any tangent vector to $M^+$ at~$p$,
then its isometric image $\tilde X = \tilde p\,JX$ at~$\tilde p$
is obtained through a rotation~$J$ by~$\sfrac\pi2$ within~$T_pM$,
followed by left translation to~$\tilde p$.
When $X$ is the (horizontal) unit tangent to~$\d M^+$ at~$p$,
the tangent plane at~$p$ is vertical and so 
$JX$ agrees with the constant conormal~$\k$ to~$\d M^+$. 
Thus $\tilde X=\tilde p\k$,
implying that the component of~$\d\plustilde M$ through~$\tilde p$
covers precisely the $\k$-Hopf circle $t\mapsto\tilde p(\cos t+\k\sin t)$.

The $\k$-Hopf circles are the fibers of the
\emph{$\k$-Hopf projection} $\Pik\colon\S^3\to\S^2$,
defined by $\Pik(q):=q\k\overline q$ in terms of quaternion multiplication. 
Thus each component of $\d\plustilde M$,
corresponding to the arc from some end~$i$ 
to end $i+1$, projects to a single point~$p_i$.
Collecting these points in~$\S^2$, we obtain an ordered $k$-tuple 
$$
  \Psi(M):=\Pik(\d\plustilde{M}) = (p_1,\ldots,p_k).
$$
As in \cite[Sect.~2]{gks} it follows that $\Psi(M)$ is
well-defined up to rotation in~$\SO(3)$.

Theorem~\ref{th:kksasym} says that the ends of a $k$-unduloid $M$ are
asymptotic to certain unduloids, defining 
\emph{asymptotic necksizes} $n_1,\ldots,n_k$ of~$M$.
These are determined by $\Psi(M)$:
\begin{theorem}\label{th:necksizes}
  Let $M$ be a coplanar $k$-unduloid of genus zero. 
  Its asymptotic necksizes $n_1,\ldots,n_k$
  are the distances between consecutive points
  of the spherical $k$-tuple $\Psi(M)=(p_1,\ldots,p_k)$:
  $$\dist_{\S^2}(p_{i-1},p_{i})=n_i.$$
\end{theorem}
\begin{proof}
  As in~\cite{gks} the result is immediate for $k=2$:
  the cousin of the \uhalf/ of an unduloid of necksize~$n$ is a spherical
  helicoid bounded by two $k$-Hopf circles;
  these circles Hopf-project to two points in~$\S^2$ at distance~$n$.  
  Then the proof given in \cite[Theorem 2.1]{gks}
  for the case $k=3$ carries over to the general case:
  the asymptotics of Theorem~\ref{th:kksasym} imply for any~$k$ that
  each end of the cousin~$\plustilde M$ has the same bounding great
  circle rays as its asymptotic unduloid cousin.
\end{proof}

\begin{corollary} \label{co:oddk}
  The set of $k$~necksizes of a coplanar $k$-unduloid obeys
  the inequalities satisfied by the edgelengths of a 
  spherical $k$-gon:
  for each subset $P\subset\{1,\ldots,k\}$ 
  with an odd number~$|P|$ of elements
  \begin{equation} \label{biswas}
    \sum_{i\in P} n_i \le \sum_{i\not\in P} n_i + \big(|P|-1\big)\pi .
  \end{equation}
\end{corollary}
\begin{proof}
  Note that the case of $|P|=k$ odd in~\eqref{biswas}
  is~\eqref{necksizesumodd}.  We prove this case only; 
  for the general case, see~\cite{bs}, where it is also shown
  that these inequalities are sharp.

  We proceed by induction,
  noting that for $k=3$ we recognize \eqref{necksizesumodd}
  as the perimeter inequality for spherical triangles~\cite{gks}.

  For $k\ge 5$, we draw minimizing geodesics of length~$a\ge 0$ 
  and $b\ge 0$ from~$p_k$ to~$p_{k-3}$ and~$p_{k-2}$, as in~\figr{kgon}.
%
\fig{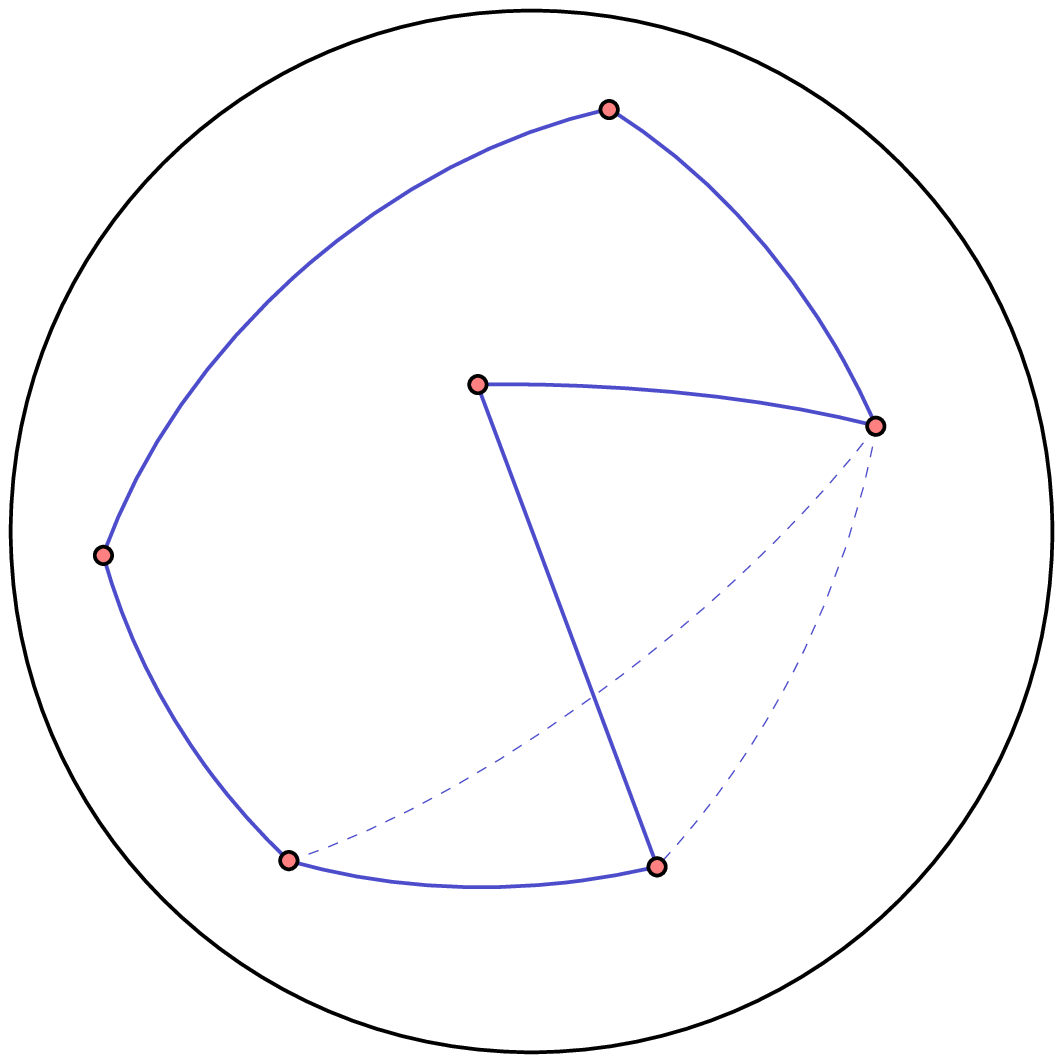}{7cm}{Perimeter inequality for spherical polygons}
{\put(83,59){$p_k$} \put(26,15){$p_{k-3}$} \put(62,15){$p_{k-2}$}
 \put(63,46){$a$} \put(77,39){$b$} \put(55,82){$p_1$}
 \put(44,14){$n_{k-2}$}
 \put(63,64){$n_k$} \put(43,66){$p_{k-1}$} \put(43,39){$n_{k-1}$} }
{Edges (dashed) of length~$a$ and~$b$ from~$p_k$
to~$p_{k-3}$ and~$p_{k-2}$ divide a $k$-gon into two triangles
and a $(k-2)$--gon; the pieces overlap in this example.
}
%
  These divide the $k$-gon into a $(k-2)$--gon and two triangles,
  which degenerate when $a$ or~$b$ vanishes.
  By induction the $(k-2)$-gon satisfies
  $$
    n_1+\cdots+n_{k-3}+a \le (k-3)\pi.
  $$
  The perimeter and triangle inequalities for the two triangles give
  $$
    n_{k-1}+n_k +b \le 2\pi
    \qquad\text{and}\qquad 
    n_{k-2}\le a+b.
  $$
  Adding these three inequalities gives~\eqref{necksizesumodd} as desired.
  Since all necksizes are positive at most $k-2$ necksizes can be~$\pi$.
\end{proof}

\begin{remark}
  This corollary extends to a broader class of 
  immersed (and not Alexandrov embedded) \cmc surfaces,
  including the surfaces constructed in~\cite{kss}.
  If $M$ has genus~$0$ with $k$~ends,
  each of which is asymptotically unduloidal, and $M$ has a mirror plane,
  then it can be shown that the disk~$M^+$ is still bounded by $k$~curves
  of planar reflection.  
  Thus the map~$\Psi$ is still defined, so that~\eqref{necksizesumodd} 
  and~\eqref{biswas} hold.
  However, inequality~\eqref{necksizesumeven} 
  depends on~\prop{dk2tk}, and need not extend to the broader class.
\end{remark}

\subsection{Spherical metrics and the classifying map}

For $k\le3$, the spherical $k$-tuple $\Psi(M)$ determines
the coplanar $k$-unduloid~$M$ uniquely. 
Since for $k>3$ this is no longer the case, we are led to introduce
a family of metrics to serve as our classifying space.

A metric is \emph{spherical} if it is locally isometric to~$\S^2$.
We will consider spherical metrics~$D$ on the oriented open two-disk.
We denote the completion of such a metric by~$\hat D$,
and refer to $\d D:=\hat D \setm D$
as the \emph{completion boundary} of $D$.
To the spherical metric $D$ we can associate 
an oriented \emph{developing map} $\phi\colon D\to\S^2$, 
unique up to rotation, which extends to~$\hat D$.
\begin{definition}
  Consider a spherical metric on the disk~$D$ 
  whose completion boundary consists of 
  fixed points $\{q_1,\ldots,q_k\}$, called \emph{vertices},
  in increasing cyclic order.
  Two such metrics are equivalent if they are related by an isometry
  extending to fix each vertex $q_1,\ldots,q_k$.
  A \emph{$k$-point metric} is an equivalence class of such metrics,
  and we let~$\Dk$ denote the space of all $k$-point metrics.
  The developing map defines a topology on~$\Dk$:
  a sequence of metrics in~$\Dk$ is said to converge if there are 
  developing maps (from the disk to~$\S^2$) which converge
  in the compact-open topology.
\end{definition}

As we noted above, if $M$ is a coplanar $k$-unduloid,
$\plustilde M$~is transverse to the $\k$-Hopf circles.
Hence $\Pik\colon \plustilde M\to\S^2$ is an immersion 
well-defined up to rotation. 
We use $\Pik(\plustilde M)$ to denote the spherical metric
on $M^+$ pulled back from~$\S^2$.
As in \cite[Lem.~3.4]{gks}, 
the exponential asymptotics of the ends of~$M$ to unduloids shows that
this spherical metric is a $k$-point metric:
The $k$~curves in~$\d M^+$ correspond to the vertices,
which develop to the $k$~labelled points~$(p_1,\ldots,p_k)=\Psi(M)$.
Thus we have shown:
\begin{theorem} \label{th:PhimapsMtoD}
  For each $k\ge2$, there is a map $\Phi\colon\Mk\to\Dk$ defined by
  $$
    \Phi(M):=\Pik(\plustilde{M}).
  $$ 
\end{theorem}
\noindent
We call this map $\Phi$ the \emph{classifying map}.
Our main theorem asserts that $\Phi$ is a homeomorphism.

\section{Decompositions of a $k$-point metric}

We analyze a $k$-point metric $D\in\Dk$ by decomposing it
into standard pieces.  
We will use this combinatorial analysis 
later to prove that the space~$\Dk$ is connected.
We assume $k\ge 2$ throughout this section.

\subsection{Polygonal metrics} \label{se:polmet}

It will simplify the exposition to consider a further
class of spherical metrics on the disk, whose completions are compact.  
\begin{definition}
  A \emph{polygonal metric} or \emph{(spherical) $k$-gon} 
  is an isometry class of spherical metrics on the open disk
  whose completion boundary is topologically a circle
  consisting of $k$~\emph{vertices} and $k$~\emph{edges} each
  developing to a great-circle arc in~$\S^2$ of length at most~$\pi$.
  We denote the space of all $k$-gons by~$\Pk$.
\end{definition}
\noindent
At each vertex~$v$ of a polygonal metric, there is an 
interior \emph{angle}~$\alpha(v)$, well-defined as a positive real number.
We introduce some further terminology:

A \emph{diagonal} in a spherical metric~$D$ is a geodesic of length
at most~$\pi$ in~$D$, whose completion in~$\hat D$ joins two vertices
(which may be consecutive).
A diagonal develops onto a minimizing great circle arc.

Two kinds of $2$-gons (see \figr{slit-lune}) will be important throughout:
A \emph{slit sphere} is 
the complement of a minimizing great circle arc in~$\S^2$.
Each of its vertices has angle $\alpha=2\pi$;
its \emph{slit length} is the length of each edge.
A \emph{lune} of angle~$\alpha$ is a $2$-gon whose vertices
develop to antipodal points of~$\S^2$; 
it has a continuous family of diagonals.

\figtwo{slit}{lune}{5cm}{Slit sphere and lune}{}{}
{Left: The complement of the great circle arc is a slit sphere.
Right: Two geodesics of length~$\pi$,
connecting a pair of antipodal points, cut the sphere into
two complementary lunes, one convex and one concave.}


A \emph{ray of spheres} is a spherical metric~$D$ whose completion boundary
is a single arc, developing to a minimizing geodesic~$\gamma$.
The preimage of~$\gamma$ decomposes~$D$ into 
an infinite union of slit spheres. 
Let us indicate an argument for this fact which is more
direct than \cite[Lemma~3.2]{gks}.  
Note that each point of a slit sphere is the endpoint
of a geodesic emanating from one of the endpoints of the slit;
each of these geodesics is minimizing and disjoint from the bounding slit.   
Thus the completeness of $D\cup\gamma$ gives
that $D$ contains a slit sphere 
attached across the completion boundary~$\gamma$.  
Continuing, by induction, we see 
$D$ contains a ray $R$ of slit spheres attached across~$\gamma$.
But $R\cup\gamma$ is complete, and thus $R$ is relatively open and closed 
in~$D$.  This shows $D=R$.

The analysis of $k$-point metrics can be reduced to $k$-gon metrics:
\begin{proposition} \label{pr:metrics}
  Given any $k$-gon, we 
  can attach rays of spheres along its edges to give a $k$-point metric.
  Conversely, any $k$-point metric $D$ can be decomposed
  into a $k$-gon~$D'$ together with $k$ rays of spheres;
  we call $D'$ a \emph{truncation} of~$D$.
\end{proposition}
\begin{proof}
  Given a $k$-gon~$D$, we can attach to each edge~$e$
  a ray of spheres with slit length equal to the length of~$e$.
  The resulting metric has no completion boundary except for
  the~$k$ vertices of~$D$, so it is a $k$-point metric.

  Conversely, given a $k$-point metric~$D$, note that its
  vertices are cyclicly ordered.  If we fix one vertex~$v$, the others
  inherit a linear ordering; let~$w$ be the last vertex in this ordering.
  We claim there is a diagonal from~$v$ to~$w$.

  First consider some reference geodesic arc~$c_0$ from~$v$ into~$D$.
  By completeness, for each angle~$\theta\in\R$ there is a
  geodesic arc~$c_\theta$ making a signed angle~$\theta$ with~$c_0$ at~$v$;
  we continue each arc~$c_\theta$ only until either
  it hits $\d D$ (at another vertex)
  or it reaches length~$\pi$ staying inside~$D$.
  Whenever the first case happens, we have found a diagonal from~$v$.

  Suppose for all $\theta$ in some interval $(\alpha,\beta)$ the second case
  occurs.  All arcs~$c_\theta$ then end inside~$D$ at the same point~$v'$
  (developing antipodally to~$v$).  Since the angle around~$v'$ is only~$2\pi$, 
  we must have $\beta-\alpha\le2\pi$,
  meaning there is never an angular gap of more than~$2\pi$ between diagonals.
  Thus we find a sequence $(d_i)_{i\in\Z}$ of diagonals from~$v$,
  indexed in increasing order of~$\theta$.  
  (If some diagonal has length~$\pi$, it is part of a continuous family
  of diagonals for some range of~$\theta$; in this case we include only
  a discrete subset of these in our sequence.)

  Let $v_i$ be the endpoint of~$d_i$.
  Since each $d_i$ cuts~$D$ into two subdiscs,
  the sequence $(v_i)$ of vertices is increasing with respect to their
  linear order.
  Thus for~$i\ge I$, the point $v_i$ is a constant vertex.
  For each~$i>I$, the loop $d_i\cup d_I$ bounds a $2$-gon~$D_i$.
  The union $\bigcup D_i =: R$ has no boundary other than~$d_I$,
  and so is a ray of slit spheres.
  Since such a ray has no further
  completion points, $v_I$ must be~$w$, the neighbor of~$v$.
  This proves the claim: each vertex has a diagonal to the consecutive vertex.

%

  Consider the set of $k$ diagonals in~$D$, connecting the pairs of
  consecutive vertices.  Since these are minimizing geodesics
  between consecutive vertices, they cannot intersect one another.
  Together, they bound a $k$-gon
  whose complement consists of $k$ rays of spheres.
\end{proof}

\begin{example}
  We exhibit polygonal metrics with
  $k-2$ edges of length~$\pi$, for each $k\ge 3$. 
  By Proposition~\ref{pr:metrics} these examples can be 
  extended with $k$~rays of spheres to form $k$-point metrics
  with $k-2$ edgelengths~$\pi$.
  By the Main Theorem, these metrics establish
  the existence of $k$-unduloids with $k-2$ cylindrical ends.

\figtwo{4gon}{4und2cyl}{5cm}{A 4-unduloid with two cylindrical ends}{
\put(26,54){$\pi$}
\put(33,31){$n$}
\put(73,46){$\pi-n$}
}{}
{Left: The solid arcs have lengths $\pi$ and $n$, while the dashed arc
has length $\pi-n$.  The $3$-gon $D_n$ develops to the upper hemisphere,
while the $4$-gon described in the text (with edgelengths $n$, $\pi$, $\pi$, $n$)
develops to the complement of the solid arcs.
Right: The classifying map takes this $4$-unduloid (with
necksizes $n$, $\pi$, $\pi$, $n$) to the $4$-point metric one
of whose truncations is the $4$-gon shown on the left. 
}
  
  Let $D_n$ be a $3$-gon developing to a hemisphere, 
  with edgelengths~$\pi$, $n$, $\pi-n$ where $0<n<\pi$.
  Reflecting across the edge of length $\pi-n$, we obtain
  a $4$-gon with edgelengths~$n$, $\pi$, $\pi$ and~$n$, as in \figr{4gon-4und2cyl}.

  We can reflect a copy of~$D_n$ across an edge of length~$n$
  of the $4$-gon to obtain a $5$-gon with
  edgelengths~$n$, $\pi$, $\pi$, $\pi$, and~$\pi-n$.
  Continuing the same way, we obtain
  $k$-gons with $k-2$ edgelengths~$\pi$; the remaining edgelengths
  are $n$, $n$ when $k$ is even, and $n$, $\pi-n$ when $k$ is odd.
  Note that the $k$~vertices of these gluings always just develop
  to the three original vertices of~$D_n$.
\end{example}

\subsection{Polygonal metrics without diagonal}

We will use diagonals to decompose a $k$-point metric in~$\Dk$
into embedded pieces.  
Let us introduce the particular $k$-gons which arise in the decomposition,
beginning with triangles ($3$-gons).
\begin{definition}
  A triangle is \emph{small} if its completion develops
  into an open hemisphere.  
  Equivalently, its developing image has area less than~$2\pi$, 
  or it has angles less than~$\pi$.  
  Its edge lengths are also less than~$\pi$.
\end{definition}

In the case of a flat metric with piecewise geodesic boundary there is a
set of diagonals which decomposes the metric into flat triangles.
This fact was used by Cos{\'\dotlessi}n and Ros 
to analyze coplanar minimal $k$-noids~\cite{cr}.
For spherical metrics, the situation is more involved.
Besides slit spheres and small triangles, our decomposition
will include pieces of the following types:
\begin{definition}
  A $k$-gon in~$\Pk$ is \emph{concave} if all its
  angles are at least~$\pi$ and no edge has length~$\pi$.
  It is \emph{embedded} if its developing
  map is an embedding into~$\S^2$ (its completion need not be embedded).
  A \emph{marked lune} is a lune of angle at most~$2\pi$
  with $k-2$ additional vertices marked on one of its edges, 
  as in \figr{mlune-longdiag};
  when the angle is less than $\pi$ we call it \emph{convex}.
\end{definition}
\figtwo{mlune}{longdiag}{5cm}{Slit sphere and lune}{}
{ \put(29,46){$p$} \put(72,46){$q$} \put(85,60){$D$} }
{Left: A lune with $k-2$ marked points on one of its edges
forms an embedded $k$-gon called a marked lune.
Right: The complement of the small convex pentagon shown on the front of the 
sphere is an embedded concave $5$-gon~$D$.  It has no diagonals.  
The great circle arc $pq$ shown, for instance, is longer than~$\pi$
and so by definition it is not a diagonal.  (A minimum-length path
from~$p$ to~$q$ exists in~$\hat D$, but it stays within $\d D$.)}

Slit spheres of slit length less than~$\pi$ are concave $2$-gons.
A triangle without diagonals is embedded, 
and is either concave or a marked lune or small.
For larger~$k$, \figr{mlune-longdiag} shows an example of 
an embedded concave $k$-gon, with no diagonals.
Another concave example is a metric isometric to a hemisphere of~$\S^2$
with $k$~vertices on its bounding equator;
it admits a diagonal if and only if some
pair of vertices is antipodal.
For slit length less than~$\pi$, slit spheres 
or chains of slit spheres are always concave.

Our goal (in \prop{conemb} below)
is to show that a concave $k$-gon with no diagonals is embedded;
for the special case of triangles this was shown in~\cite{gks}.
To formulate a lemma, let us fix a vertex~$v$ of~$D$
and consider geodesics in~$D$ emanating from~$v$,
as in the proof of Proposition~\ref{pr:metrics}.
For each angle~$\theta\in\big(0,\alpha(v)\big)$,
there is a unique geodesic arc~$c_\theta$ making an angle~$\theta$ 
with the edge after~$v$; again we continue~$c_\theta$ only until
either it hits $\d D$, or it reaches length~$\pi$ staying inside~$D$.
Thus the only points of~$\d D$ contained 
in an arc $c_\theta$ are its initial point~$v$ and perhaps its endpoint.

\begin{lemma}\label{lem:twocase}
  Let $D$ be a $k$-gon which is not isometric to a hemisphere, and
  suppose that $D$ has no diagonals. 
  Then all edges of $D$ have length less than~$\pi$.
  Moreover, at each vertex $v$ the following alternative holds.  Either
\begin{itemize}
\item[\ia] all arcs~$c_\theta$ can be continued to length~$\pi$, 
  where they converge at a point $v'\in D$
  which develops antipodally to~$v$, or 
\item[\ii] all arcs~$c_\theta$ have length at most~$\pi$ 
  and end on the interior of the same edge~$e$ (not incident to~$v$),
  where they meet~$e$ transversely.
\end{itemize}
  The second case can only arise when $\alpha(v)<\pi$. 
\end{lemma}

\begin{proof}
  If some edge $e$ has length~$\pi$, 
  then there is a sufficiently thin lune contained in~$D$ 
  whose boundary includes~$e$.
  But such a lune contains a continuous family of diagonals, 
  joining the endpoints of~$e$.  
  This contradiction shows the edges have length less than~$\pi$.
  
  Let us now prove the alternative.
  Since~$D$ has no diagonals, no arc~$c_\theta$ can hit~$\d D$ at a vertex.
  For any fixed~$\theta$,
  either~$c_\theta$ reaches length~$\pi$ without encountering~$\d D$,
  or $c_\theta$ first hits~$\d D$ in an edge $e$ transversely.
  Since the edges have length less than~$\pi$, the edge~$e$
  in the second case is one of the $k-2$ edges not incident to~$v$.
  We claim that each of these $k-1$ possibilities occurs for
  an open set of angles~$\theta$.  
  Since $(0,\alpha(v))$ is connected, it then follows that we must have
  the same endpoint behavior for all angles~$\theta$.  

  To prove the claim, consider~$c_\theta$.  Let~$X$ be~$\d D$ with~$v$
  and its two incident edges removed; if $c_\theta$ ends on
  an edge~$e$ we also remove~$e$.  Then 
  the two disjoint closed subsets~$c_\theta$ and~$X$ of~$\hat D$
  have a positive distance~$\delta$.
  Within the $\delta$-neighborhood of~$c_\theta$,
  we can construct the required nearby geodesics starting from~$v$
  and with the same endpoint behavior.

  It remains to show that for $\alpha(v)\ge\pi$, case \ii\ cannot occur.
  If $\alpha>\pi$ then we have two opposite arcs from~$v$ hitting
  the same boundary edge~$e$ transversely.  The two intersection points
  of this great circle with~$e$ must be antipodal, contradicting
  the fact that $e$ has length less than~$\pi$.
  To rule out $\alpha=\pi$, we similarly use two limit arcs in~$\hat D$
  which hit the same edge~$e$.  If~$e$ is transverse to the limit
  arcs, we get the same contradiction.  Otherwise $e$ lies along the
  limit arcs; however, by assumption $D$ is not a hemisphere, so in
  this case one of the limit arcs must include a diagonal of~$D$,
  a contradiction.
\end{proof}

\begin{proposition}\label{pr:concave}
  If a $k$-gon $D\in\Pk$ has no diagonal,
  then either $D$ is concave or it is a small triangle. 
\end{proposition}
\begin{proof}
  For $k=2$, the $2$-gons without diagonals are slit spheres,
  and hence concave.  So we assume $k\ge 3$.
  We first consider the case that all angles of~$D$ equal~$\pi$.
  The developing image of the completion boundary is then contained
  in a single great circle and, by the Gauss--Bonnet formula, 
  the metric must have area~$2\pi$.  
  Therefore $D$ is isometric to a hemi\-sphere.
  It cannot have an edge of length~$\pi$ for then
  it would have diagonals near this edge;
  thus it is concave and we are done.

  We may thus assume that $D$ is not a hemisphere,
  meaning the hypotheses of Lemma~\ref{lem:twocase} hold.
  We suppose $D$ is not concave and claim $D$ must be a small triangle. 
  In the following we use $\overline{vw}$ to denote
  the edge between consecutive vertices~$v$ and~$w$, 
  not including its endpoints.

  There is at least one vertex~$w$ with $\alpha(w)<\pi$.
  Consider the consecutive vertices $t,u,v,w,x$ 
  (not necessarily distinct for small~$k$). 
  We apply the lemma to the arcs~$c_\theta$ at~$v$. 
  Because the length of~$\overline{vw}$ and the angle at~$w$ are
  both less than~$\pi$, as in \figr{nodiag},
  arcs~$c_\theta$ near~$\overline{vw}$ end on edge~$\overline{wx}$.
  The lemma then says all $c_\theta$ end on~$\overline{wx}$;
  this implies that we also have $\alpha(v)<\pi$.
\fig{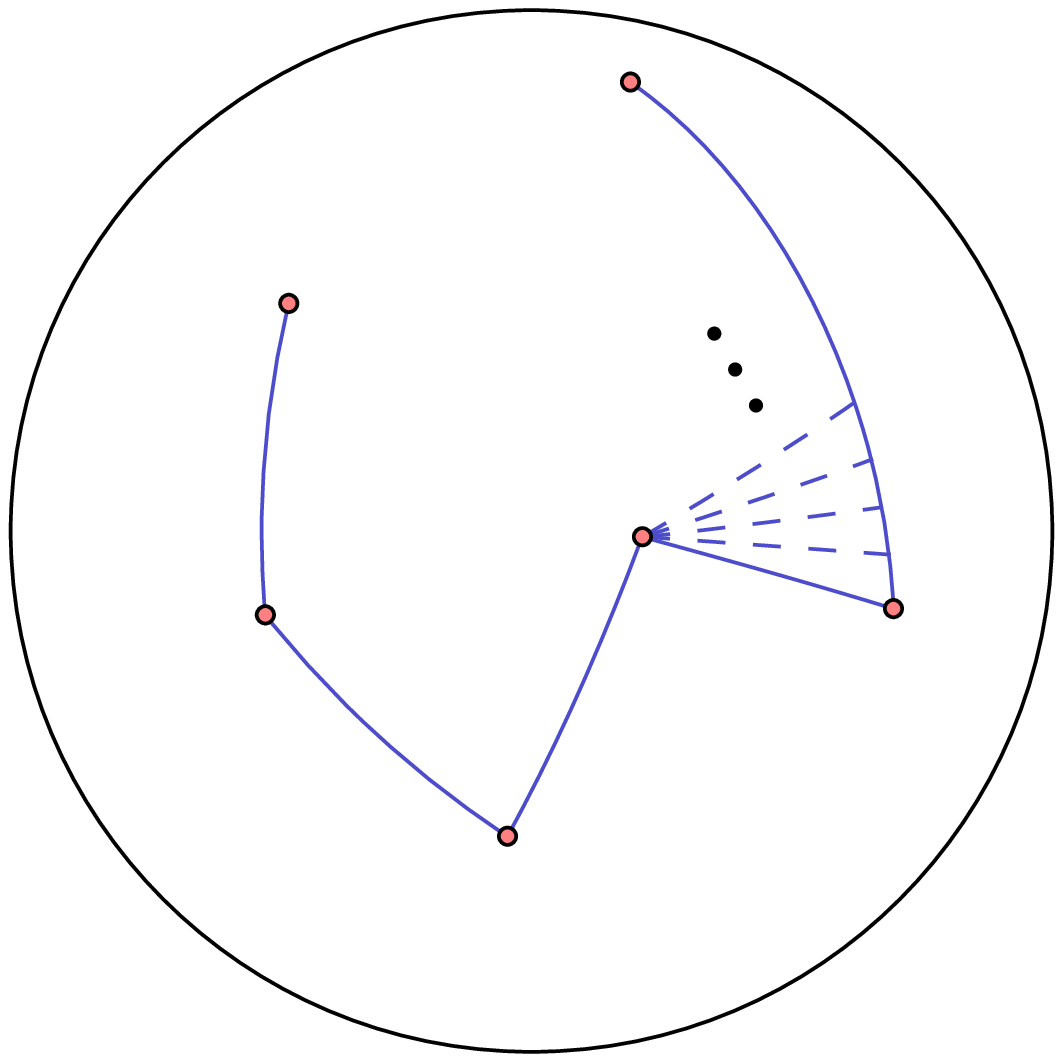}{6cm}{Polygonal metrics with no diagonals}
{ \put(21,40){$t$} \put(46,17){$u$}
  \put(56,50){$v$} \put(81,38){$w$} \put(54,90){$x$} }
{The rays $c_\theta$ from $v$ initially hit~$\overline{wx}$
so by Lemma~\ref{lem:twocase} they all do, and we deduce
that, unlike in this picture, $\alpha(v)<\pi$.}
  
  Repeating this argument from~$v$, we find next that $\alpha(u)<\pi$.
  Now consider again the arcs $c_\theta$ at~$v$.
  They have length less than~$\pi$ and
  end on the edge~$\overline{wx}$ when~$\theta$ is small, while they 
  end on the edge~$\overline{tu}$ for~$\theta$ close to $\alpha(v)$.
  By the lemma, we must have $\overline{wx}=\overline{tu}$, 
  implying that $k=3$; moreover all edges and angles are less than~$\pi$.  
  That is, $D$ is a small triangle, as claimed.
\end{proof}

We can characterize concave $k$-gons without diagonals further:
\begin{proposition} \label{pr:conemb}
  A concave $k$-gon with no diagonals is embedded.  
\end{proposition}
\begin{proof}
  It follows as in the previous proof that 
  if all vertex angles equal~$\pi$ then
  equivalently the metric is a hemisphere and hence embedded.
  Otherwise Lemma~\ref{lem:twocase} is applicable. 

  We first claim that $\alpha(w)\in[\pi,2\pi]$ for each vertex~$w$.
  Since $D$ is concave, $\alpha(w)\ge\pi$,
  so if we apply Lemma~\ref{lem:twocase} we are in case~\ia.
  That is, all arcs~$c_\theta$ from $w$ continue to length~$\pi$,
  where they end at a common point $w'\in D$.  We call their
  union the \emph{vertex lune} of~$D$ at~$w$.  Since $w'$ is an
  interior point of~$D$, and~$D$ immerses into~$\S^2$ under the developing map, 
  as in the proof of Proposition~\ref{pr:metrics},
  the angle $\alpha(w)$ of this lune cannot exceed~$2\pi$.
  This proves the claim.

  We then use induction on~$k$.
  For $k=2$, we know a $2$-gon without diagonals is an embedded slit sphere.  
  For the induction step, we will show that given a $k$-gon~$D$,
  we can join a triangle to it to obtain a $(k-1)$--gon~$D'$ that still
  has no diagonals.  By the induction hypothesis, $D'$ is embedded, 
  and so the same must hold for its subset~$D$.

  Since the $k$-gon is not a hemisphere
  we can choose a vertex~$v$ with $\alpha(v)\in(\pi,2\pi]$.
  Let us assume for now that $\alpha(v)<2\pi$.
  Then the three consecutive vertices $u,v,w$ of~$D$
  develop to the three vertices of a small
  triangle~$\Delta$, which is locally exterior to the image of~$D$.
  We define a polygonal metric~$D'\in\D_{k-1}$, as in \figr{concind},
  isometric to the union of~$D$ with~$\Delta$, glued across
  the two edges incident to~$v$.
  More precisely, we glue
  the closure of~$\Delta$ to $D$ and then remove~$\overline{uw}$.
  The completion boundary of~$D'$ agrees with that of~$D$, 
  except that $\overline{uv}$ and~$\overline{vw}$ are
  replaced by~$\overline{uw}$.
\fig{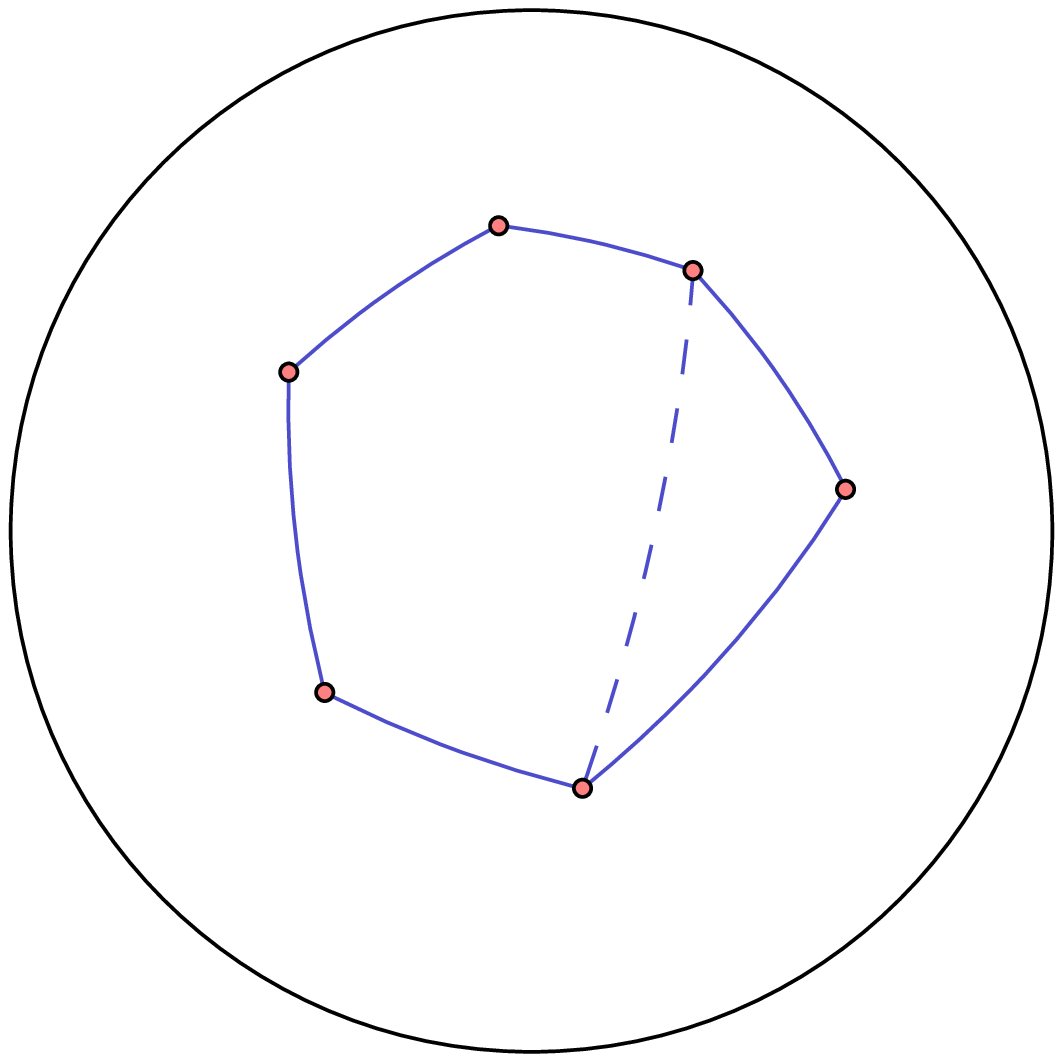}{6cm}{Concave metrics with no diagonal are embedded.}
{ \put(57,23){$u$} \put(81,50){$v$} \put(64,77){$w$}
  \put(77,28){$D$} \put(67,50){$\Delta$} }
{By induction, the $(k-1)$--gon~$D'$, obtained by gluing
the triangle~$\Delta$ to the concave $k$-gon $D$, is embedded.}

  The angles at~$u$ and~$w$ have only increased under this gluing,
  so $D'$ is still concave.
  We claim that $D'$ does not have a diagonal.
  To see this, consider the open hemisphere $H\subset D'$,
  bounded by $\overline{uw}$, and contained 
  in the union of~$\Delta$ with the vertex lune at $v\in D$.
  No vertices of $D'$ are in~$H$, so a
  diagonal of~$D'$, being a minimizing geodesic with endpoints outside~$H$,
  stays outside~$H$.  Hence any diagonal of~$D'$ stays outside~$\Delta$
  and is thus also a diagonal of~$D$.
  However, by assumption there are no such diagonals.
  This completes the inductive proof under the assumption that $\alpha(v)<2\pi$.

  When $\alpha(v)=2\pi$, essentially the same ideas apply, except that
  $\Delta$ degenerates to a segment.  When $u$ and $w$ develop to
  distinct points, we obtain exactly as above a metric~$D'\in\P_{k-1}$
  with no diagonals.  When~$u$ and~$w$ develop to the same point,
  note that we must have $k\ge4$.  Here we get $D'\in\P_{k-2}$
  by gluing the edge $\overline{vw}$ (being the degenerate triangle~$\Delta$) 
  to~$D$, and removing the vertex $u=w$.  Again by the same argument,
  the metric~$D'$ cannot have any diagonals.
\end{proof}

Note that an embedded concave polygon has diagonals only if it has
two antipodal vertices; in this case it is a lune with angle in $[\pi,2\pi]$,
possibly with extra vertices marked along its edges.

\subsection{Decomposition of polygonal and $k$-point metrics}

We will now use the previous proposition to decompose any polygonal
metric into a finite union of embedded pieces.
Each piece is either a small triangle,
or an embedded concave polygon with no diagonals,
or a marked lune of angle less than~$\pi$.
It will follow immediately that any
$k$-point metric can be decomposed into the same pieces.

\begin{theorem}  \label{thm:decomp}
  Any polygonal metric $D\in\Pk$ can be decomposed
  by diagonals into a finite union of slit spheres,
  small triangles, embedded concave polygons without diagonals,
  and marked lunes.
\end{theorem}
\begin{proof}
  By Propositions~\ref{pr:concave} and \ref{pr:conemb}, 
  if $D$ has no diagonals, it is already one of the allowed pieces.  
  Otherwise, we work by induction on~$k$, using diagonals to simplify $D$.  
  Suppose that $D$ has a diagonal~$d$ connecting two nonconsecutive vertices. 
  In this case the induction step is immediate
  since $d$ splits $D$ into two pieces each with fewer edges. 

  The remaining case is when $D$ has a diagonal connecting consecutive
  vertices, but no diagonals connecting nonconsecutive ones.
  This includes the base case ($k=2,3$) of the induction on~$k$.
  Here we will establish the decomposition directly, without using 
  the induction hypothesis.
  Instead we work by induction on the area of~$D$ (in multiples of~$4\pi$).

  So we have a metric~$D$, of area less than $4\pi n$ for some~$n\in\N$,
  where all diagonals connect consecutive vertices.
  Choose a shortest diagonal~$d$, and call its endpoints $v$ and~$w$.
  Suppose first~$d$ has length less than~$\pi$.
  Then the $2$-gon bounded by~$d$ and the edge $\overline{vw}$ of~$D$
  is a union of $m$ slit spheres, where $1\le m\le n$.  
  The $k$-gon on the other side of~$d$ has area $4\pi m$ less than 
  that of~$D$.  Thus for the base case $n=1$ this is impossible, 
  meaning that no such diagonal can exist so that $D$ itself is embedded.
  On the other hand, for $n>1$ the metric decomposes into the allowed 
  embedded pieces by induction.
  
  The more difficult case is when~$d$ has length~$\pi$.
  In this case we consider the lune~$L$ of maximal angle which is
  contained in~$D$ and incident to the edge~$\overline{vw}$.
  Then the other edge~$d'$ bounding~$L$ must include some point of~$\d D$.
  Since there are no diagonals shorter than~$\pi$,
  the arc~$d'$ must lie entirely in~$\d D$.
  Thus $D$ is a marked lune, though with angle~$\alpha$
  possibly greater than~$2\pi$.
  It has a continuous family of diagonals,
  so we can use these to decompose $D$ into $n$~slit spheres
  and one marked lune of angle at most~$2\pi$ (adjacent to~$d'$).
  Thus also in this case we have arrived at a decomposition into
  the desired embedded pieces.
\end{proof}

Given a $k$-point metric $D$ we can 
remove its rays as in Proposition~\ref{pr:metrics} 
and then apply Theorem~\ref{thm:decomp} to obtain the following 
decomposition:
\begin{corollary}\label{co:kptdec}
  Any $k$-point metric $D\in\Dk$ can be decomposed into $k$ rays
  of slit spheres together with a finite union of 
  slit spheres, small triangles, embedded concave polygons
  with no diagonals, and marked lunes. 
\end{corollary}
When applied to the Hopf projection $D=\Phi(M)$ of
a coplanar $k$-unduloid $M\in\Mk$ this decomposition 
can be regarded as a decomposition of~$M$ into bubbles.

There is an equivalent reformulation of the decomposition in terms of
graphs dual to the tessellating pieces.  This graph is a tree,
topologically equivalent to the skeletal graphs considered earlier by
Schoen~\cite{sch}, Kusner~\cite{kmsri}, and Kapouleas~\cite{kap}.

Note that the decomposition of a $k$-point metric in general is not unique. 
Clearly, the location of diagonals of length $\pi$ is not determined;
in addition, many metrics admit combinatorially different choices of
diagonals, with the simplest example being a convex square.  

In Proposition~\ref{pr:metrics} we removed rays (of slit spheres)
from a $k$-point metric.
In view of the decomposition established by the corollary above,
we can make this more precise:
\begin{definition}
  Given a $k$-point metric $D\in\Dk$, a \emph{small truncation}
  is a metric $\Dbox\in\Pk$ obtained from~$D$
  by removing (closed) rays so that no slit 
  sphere remains incident to any bounding edge.
\end{definition}
When $D$ has no pair of consecutive vertices at distance~$\pi$, 
the small truncation is uniquely defined.  In the general case,
the ambiguity is by less than one slit sphere at each of the $k$~ends.  

\figr{ex6pt} illustrates the results of this section.
\fig{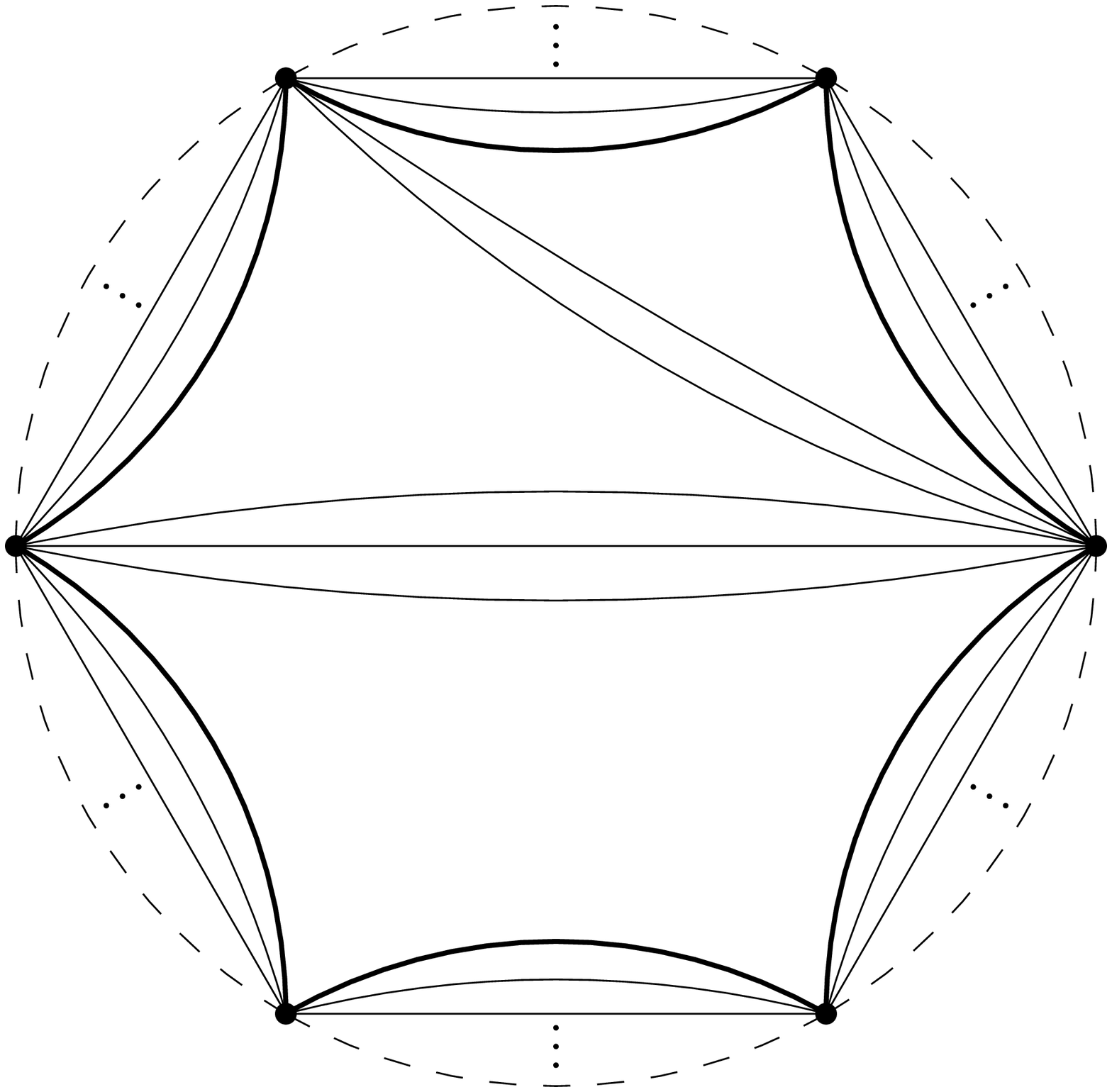}{7cm}{A six-point metric on the open disk}
{}{This $6$-point metric~$D$ on the open disk is decomposed
    by the six thick arcs into its small truncation $\Dbox$
    and six rays of slit spheres.  In turn, the polygon
    metric $\Dbox$ is, as in \cor{kptdec}, decomposed into two
    triangles, three slit spheres, and one $4$-gon (necessarily
    embedded, concave, and without diagonals).  A $6$-point
    metric with such a decomposition will correspond (under
    our classifying map~$\Phi$) to a $6$-unduloid with two
    midsegments, one having two unduloid bubbles and the
    other having just one.}

\section{The space of $k$-point metrics}

We now show that $\Dk$ is a connected manifold. 
As for $k=3$ we consider the projection~$\Pi$ from $\Dk$ to $k$-tuples 
of points on $\S^2$, given by restricting the developing map 
of a $k$-point metric to its vertices.
For $k=3$ the mapping $\Pi$ is a diffeomorphism from~$\D_3$ 
to the open 3-ball of triples in~$\S^2$, as in~\cite{gks}.
For larger~$k$ this is only true locally.  
From now on we assume $k\ge 3$.

\subsection{$\Dk$ is a manifold}

Let us define ${\mathcal S}_k$ to be the set of $k$-tuples
of sequentially distinct points in~$\S^2$:
$$
  {\mathcal S}_k:=\{(p_{1},p_{2},\ldots,p_{k})\in (\S^2)^k, 
  p_{i}\not=p_{i+1} \text{ for }i=1,\ldots,k\}.
$$
Define ${\mathcal A}_k\subset {\mathcal S}_k$ to be the set 
$\{(p,q,\ldots,p,q): p\not=q\in\S^2\}$ for $k$~even and the empty set
for $k$~odd.
Note that ${\mathcal A}_k$ contains all $k$-tuples in ${\mathcal S}_k$ 
consisting only of a pair of antipodal points (namely, when $q=-p$).
Thus the action of $\SO(3)$ is free on 
the manifold ${\mathcal S}_k\setm {\mathcal A}_k$, 
and the quotient
$$
  \T_k:=\big({\mathcal S}_k\setm {\mathcal A}_k\big) \big/ \SO(3),
$$
is a manifold.

\begin{proposition}\label{pr:dk2tk}
  For $k\ge3$, the completion boundary of any $k$-point metric develops
  to a $k$-tuple in $\Tk$.  
  This gives a surjective map $\Pi\colon \Dk\to\Tk$,
\end{proposition}
Note, however, that for $k=2$ the image of $\D_2$ includes
$(p,-p)\in {\mathcal A}_2={\mathcal S}_2$.
This fact explains why $\D_2\isom(0,\pi]$ fails to be a manifold
but instead has a boundary point~$\pi$.  Of course,
$\pi$~corresponds to the cylinder in $\M_2\isom\D_2$,
where rotations also fail to act freely.
\begin{proof}
  For $D\in\Dk$, Proposition~\ref{pr:metrics} gives
  a $k$-gon $\Dbox\in\Pk$, obtained by removing rays from~$D$.
  The spheres of each ray have nonzero slit length, 
  and so consecutive vertices must develop to distinct points. 
  Thus any developing image of the completion boundary $\d D$
  lies in ${\mathcal S}_k$.

  To show the image is not in ${\mathcal A}_k$, suppose
  to the contrary that the completion points of $D\in\Dk$ develop 
  onto just one pair $p\not=q$ of points in~$\S^2$.
  Then the decomposition of Corollary~\ref{co:kptdec}
  can only contain slit spheres of slit length $\dist(p,q)$.
  Such slit spheres
  cannot be joined together to produce a $k$-gon for $k>2$.

  Finally we show~$\Pi$ is surjective.  
  For $k=3$ this is obvious, while for $k>3$ this follows by induction:
  There exists a nonconsecutive pair $p_i,p_j$ of vertices
  with distinct images.  By joining two metrics along a
  diagonal developing to~$\overline{p_ip_j}$
  we obtain the desired $k$-point metric.
\end{proof}
\begin{corollary} \label{co:evenk}
  For $k$ even, at least two ends of a $k$-unduloid
  are non-cylindrical, and so $n_1+\cdots+n_k < k \pi$ holds.
\end{corollary}
\begin{proof}
  If $k-1$ ends were cylindrical,
  then the $k$-tuple~$\Psi(M)$ 
  would alternate between a pair of antipodal points
  (making also the $k^{\text{th}}$ necksize cylindrical).
  But this is impossible by~\prop{dk2tk}.
\end{proof}

The map~$\Pi$ relates the $k$-tuple map~$\Psi$ and classifying
map~$\Phi$ introduced in Section~1:  For any $k$-unduloid $M\in\Mk$,
we have $\Phi(M)\in\Dk$ and $\Psi(M)=\Pi(\Phi(M))\in\Tk$.

Let $\Tk^*\subset\Tk$ denote the open set of 
all pairwise distinct $k$-tuples in $\S^2$ (mod $\SO(3)$),
and let $\Dk^*:=\Pi^{-1}(\Tk^*)\subset\Dk$.
For $k=3$ we have $\T_3^*=\T_3$ and so $\D_3^*=\D_3$.

\begin{theorem} \label{th:mfld}
  For any~$k\ge3$, the map $\Pi\takes\Dk\to\Tk$ is a covering map over $\Tk^*$.
  It is a local diffeomorphism everywhere, implying that~$\Dk$ is a
  real-analytic manifold of dimension $2k-3$.
\end{theorem}
\begin{proof}
  Consider a $k$-point metric $D\in\Dk$, which maps to some $k$-tuple 
  $T\in\Tk$.   Pick any representative $(v_1,\ldots,v_k)\in(\S^2)^k$
  for $T$; the $v_i$ are the images of the completion boundary points $q_i$
  of~$D$ under some developing map~$\phi$.

  By introducing a finite set of additional vertices
  we can find a triangulation~$\Delta$ of~$\S^2$ by small triangles
  whose vertices include the~$v_i$.  Passing to a
  subdivision if necessary, we may assume that the~$v_i$ are not
  adjacent in~$\Delta$, in other words that their 
  star neighborhoods~$S_i\subset\S^2$ are disjoint when $v_i\ne v_j$.  
  (Note that we may have $v_i=v_j$ for $i\ne j$ if $T\not\in\Tk^*$.)
  We claim that $\Pi$ is a local diffeomorphism from a neighborhood of~$D$ 
  onto $U:=(S_1\times\cdots\times S_k)/\SO(3)$.  
  It follows that $\Dk$ is a manifold of dimension $2k-3$; 
  we endow it with the real-analytic structure coming from these charts.

  To prove the claim, note that the pullback $\phi^{-1}(\Delta)$
  gives a triangulation of the metric space~$\hat D$.  For any
  $(w_1,\ldots,w_k)\in S_1\times\cdots\times S_k$, we construct
  a completed $k$-point metric with that developing image as follows.
  Each triangle of $\phi^{-1}(\Delta)$ incident to the completion
  boundary~$\d D$ has exactly one vertex $q_i\in\d D$.  
  We continuously change the metric
  on any such triangle according to the motion of $v_i$ to~$w_i$, but we
  fix the metric on the remaining triangles (those not incident to $\d D$).
  Under this new $k$-point metric, $q_i$ develops to~$w_i$; this confirms
  that we have constructed a continuous local inverse for~$\Pi$.

  To show that $\Pi$ is a covering over $\Tk^*$, note that $\Tk^*$ is
  exactly the set where the $v_i$ are distinct.  In this case, the
  neighborhood $U=(S_1\times\cdots\times S_k)/\SO(3)$ is trivially
  covered by $\Pi^{-1}(U)$ because we can use the same triangulation~$\Delta$
  for any $D\in\Pi^{-1}(U)$.
\end{proof}

\subsection{$\Dk$ is connected}

\begin{figure}
\hspace*{-9mm}
\begin{overpic}[width=.32\textwidth]{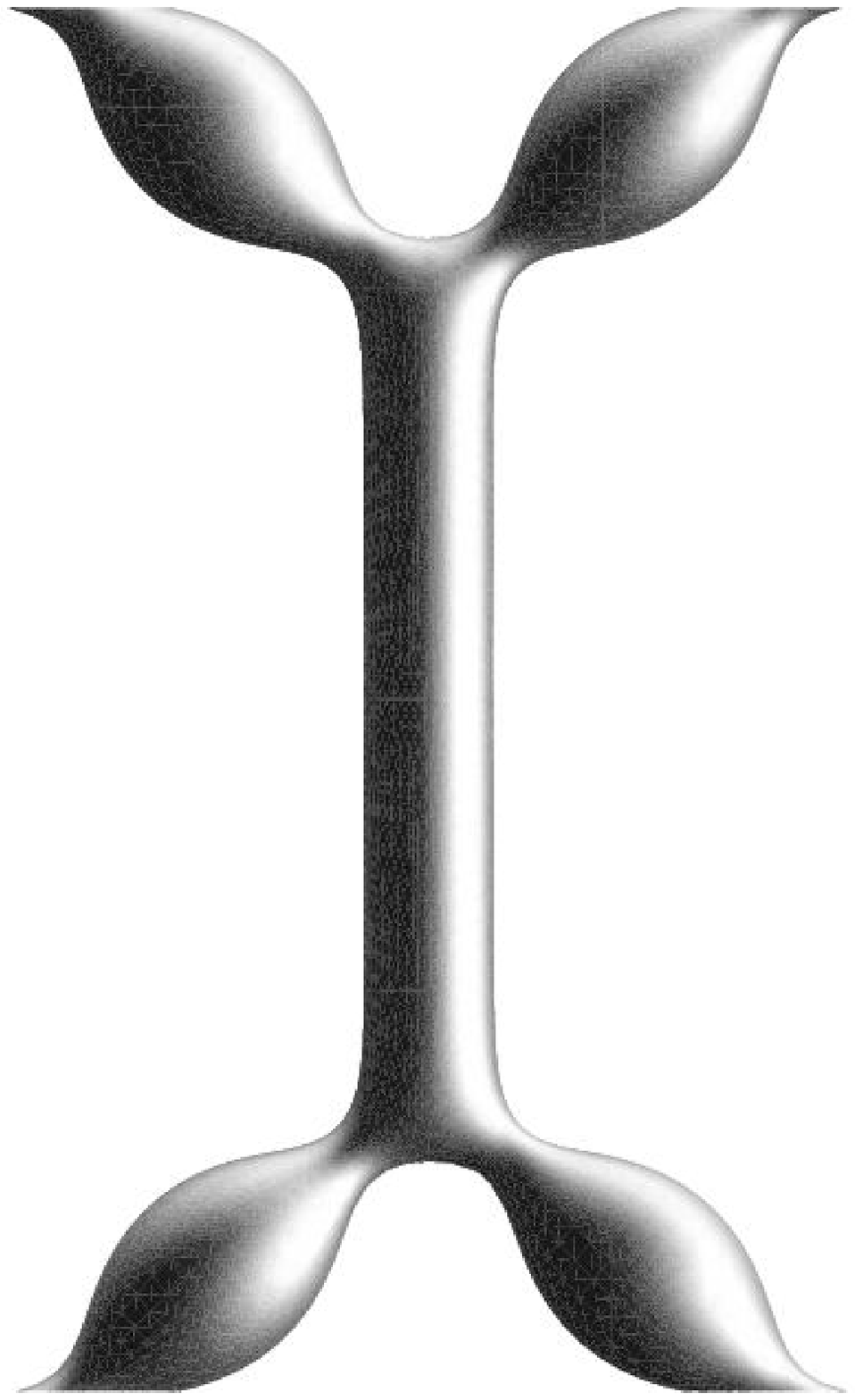} \end{overpic}
\begin{overpic}[width=.32\textwidth]{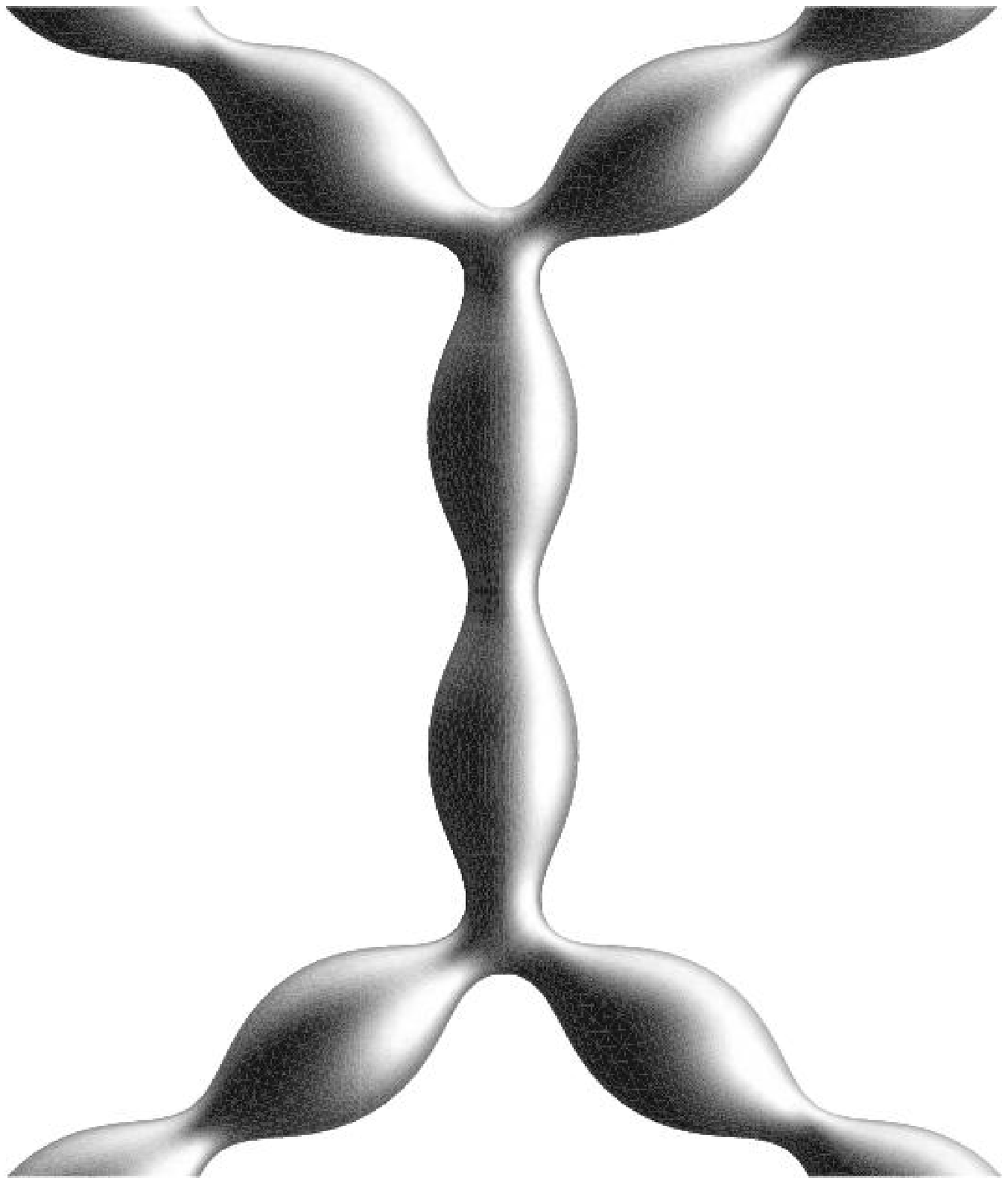} \end{overpic}
\hspace{3mm}
\begin{overpic}[width=.32\textwidth]{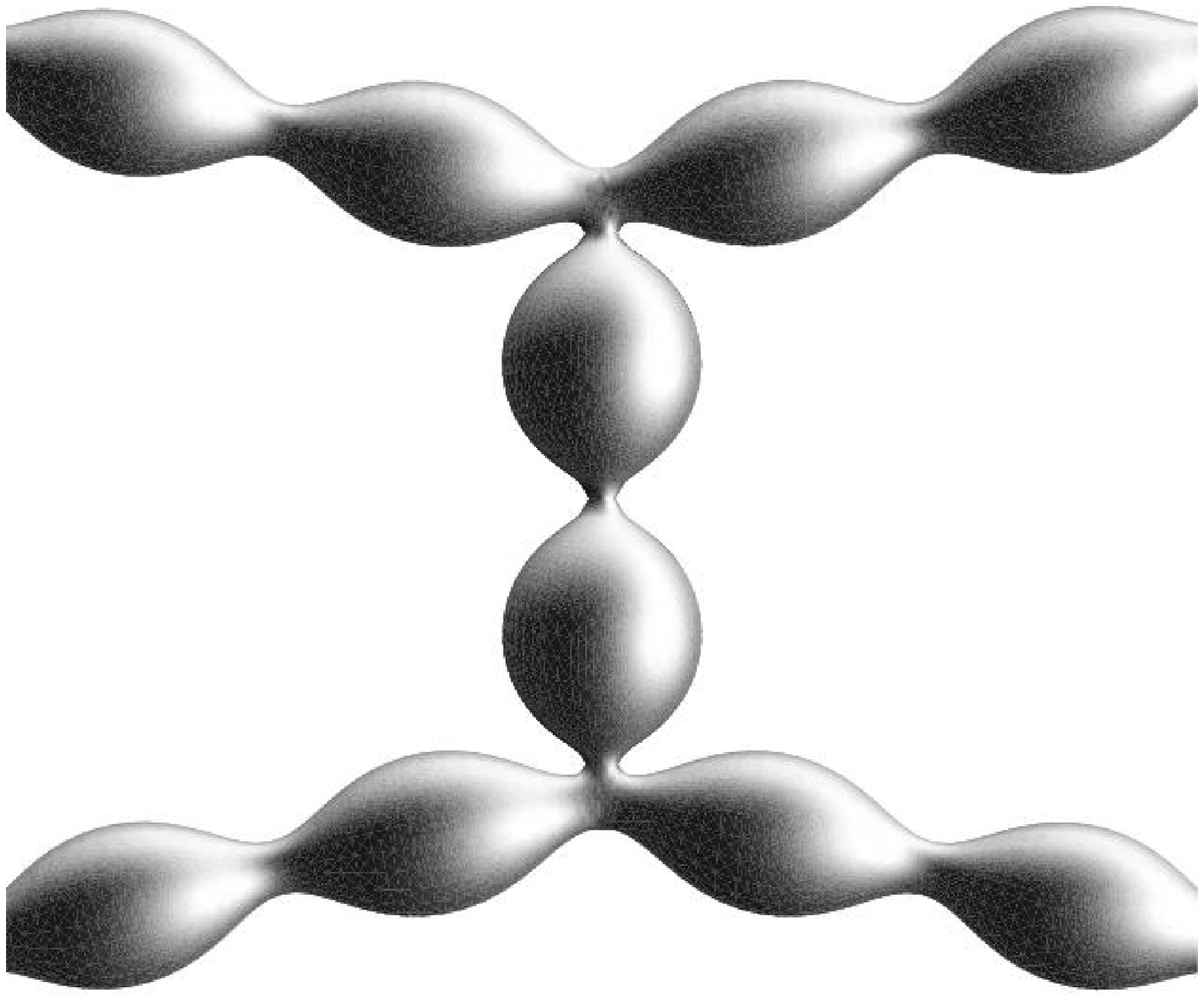} \end{overpic}
\caption[Curvature blow up]{
    This family of $4$-unduloids with fixed necksizes $\pi/2$
    exhibits a curvature blow-up on the three necks attached to the
    two central bubbles.  Its singular limit is a union of two unduloids 
    and two spheres.  Under~$\Psi$, the surfaces Hopf project
    to the two poles of~$\S^2$ as well as two points on the equator 
    which approach  one another; thus the vertex distance has limit~$0$ 
    and Lemma~\ref{le:curv} does not apply.
}\label{fig:curvature blow up}
\end{figure}

What initially led us to believe that~$\Mk$ is connected
was that we can imagine a path from any given $k$-unduloid $M\in\Mk$
to a standard one, say a dihedrally symmetric one.
In particular, we could deform~$M$ such that one of its midsegments becomes
cylindrical (as in \figr{curvature blow up} viewed from right to left),
and then (as in \figr{rhombic with cyl segment})
telescope inward along that midsegment. 
Repeating this process until all midsegments have been retracted, 
we could finally adjust the necksizes of the ends until they are equal. 

This intuitive idea will guide our proof of connectedness of~$\Dk$.
Our first lemma singles out a class of $k$-point metrics 
for which the decomposition of Corollary~\ref{co:kptdec} is unique.
Let us introduce the subset $\E_k\subset\Tk$ 
of equatorial $k$-tuples subject to the condition
that all points occur in counterclockwise order on the equator,
and consecutive pairs have distances in $(0,\pi)$.
\begin{lemma} \label{le:equatorial}
  Suppose that $D\in\Dk$ is a metric such that $\Pi(D)\in\E_k$ 
  and $D$ has no diagonal of length~$\pi$. 
  Then the decomposition of~$D$ asserted in Corollary~\ref{co:kptdec} 
  is unique.  
  More specifically, the small truncation $\Dbox$ is a union of
  some number $m(D)\ge 0$ of pieces isometric to slit spheres
  and exactly one piece isometric to a hemisphere. 
\end{lemma}
\begin{proof}
  The decomposition pieces 
  can only be slit spheres or embedded concave $k$-gons without diagonals.
  Let us consider the developing image of~$D$.
  The diagonals of~$D$ can only develop to the equator, 
  and so the concave $k$-gons are isometric either to hemispheres 
  or to slit spheres.   Moreover, all rays of slit spheres
  are joined from the same 
  side to the equator, meaning that 
  the small truncation $\Dbox$ has area an odd multiple of~$2\pi$.
 
  Thus there is at least one
  decomposition piece isometric to a hemisphere.
  The orientation of~$D$ means that any such hemisphere~$H$
  has its bounding vertices in increasing cyclic order;  
  however, the developing image of~$\d H$ covers the equator once,
  so~$H$ must develop to the upper hemisphere. 
  If there were two or more hemispheres in the decomposition,
  we could find a pair 
  connected by a chain of pieces isometric to slit spheres;
  this pair would develop to complementary hemispheres, contradicting
  the observation above. 

  
  Let us now show the decomposition is unique.
  We use preimages of the equator to subdivide~$D$ into hemispherical tiles,
  developing alternately to the upper and lower hemisphere.
  This canonical decomposition is the same as the one above, except that
  each slit sphere has been divided (along an internal arc of length
  greater than~$\pi$) into two hemispheres.
  Thus the hemisphere piece in the original decomposition is distinguished
  as the unique tile which has all edge lengths less than~$\pi$.
\end{proof}

We can now prove the following:
\begin{theorem}\label{th:connected}
  The manifold $\Dk$ is connected.
\end{theorem}
\begin{proof}
  By Theorem~\ref{th:mfld} the complement of~$\Dk^*$
  has codimension~2 in~$\Dk$; for $k=3$ it is empty. 
  Thus proving $\Dk$ is connected is equivalent 
  to showing that $\Dk^*$ is connected.
  %
  %
  To prove the latter, it is sufficient to pick some $T\in\Tk^*$ and
  show that the deck transformation group of the covering
  $\Pi\colon\Dk^*\to\Tk^*$
  acts transitively on the particular fiber~$F:=\Pi^{-1}(T)\subset\Dk^*$.

  To define~$F$, we pick a $k$-tuple $T\in\E_k$ 
  with no pairs of points antipodal.
  For instance, if $k$ is odd, 
  then $T$ could be $k$~points with maximal dihedral symmetry,
  while for even~$k$, 
  a slight perturbation of this configuration within~$\E_k$ would work.  
  Since a metric $D\in F$ cannot have diagonals of length~$\pi$,
  Lemma~\ref{le:equatorial} is applicable to~$D$ and gives a well-defined
  number~$m(D)$ of pieces isometric to slit spheres
  in the unique small truncation~$\Dbox$.
  Note there is a unique metric in~$D_0\in F$
  with decomposition number~$m(D_0)=0$, isometric to 
  a hemisphere with $k$~rays of spheres attached.

  For $k=3$, any $D\in F$ has $m(D)=0$ since the
  small truncation~$\Dbox$ is the unique triangle determined by $\Pi(D)\in\T_3$.
  This shows $\D_3^*=\T_3^*$; but $\T_3^*=\T_3$ is an open 3-ball 
  (see~\cite{gks}). 
  So let us from now on assume $k\ge 4$.

  For $D\in F$ we will construct a loop in~$\Tk^*$ based at~$T$
  whose lift in~$\Dk^*$ runs from $D$ to~$D_0$.  
  We can assume $D$ has decomposition number $m(D)>0$;
  by induction it is sufficient to construct a loop~$\gamma$
  whose lift~$\tilde\gamma$ connects~$D$ to some $D'\in F$
  for which~$m$ has decreased by one.
  We will choose~$\gamma$ as the concatenation $\alpha * \lambda * \alpha^{-1}$
  of three subpaths, where $\lambda$ is a loop in~$\Tk^*$
  and $\alpha$ is a path in~$\E_k$.
  On the level of \CMC/ surfaces, the lift of~$\lambda$ telescopes
  a midsegment, while $\alpha$ adjusts necksizes.

  Since $m(D)>0$, the metric $D$ has at least one
  diagonal $d$ between two nonconsecutive vertices $v_i$, $v_j$.
  There is a path $\alpha$ in $\E_k$ from $T=\Pi(D)$ 
  to some $k$-tuple $S\in\E_k$ in which~$v_i$ and~$v_j$ are antipodal.
  Among pairs of vertices joined by a diagonal,
  we can assume that no pair becomes antipodal along our path $\alpha$, 
  since otherwise we instead take the first such diagonal to be~$d$,
  and replace $\alpha$ by the shortened path.

  The path~$\alpha$ has a unique lift $\mtilde\alpha$ in~$\Dk^*$ through~$D$.
  We let $E\in\Pi^{-1}(S)$ be its terminal metric.
  By Lemma~\ref{le:equatorial} the decomposition of each metric
  along~$\mtilde\alpha$ is unique.
  Hence up until the terminal time, the diagonal~$d$ as well as the
  entire decomposition continues uniquely along the lift~$\mtilde\alpha$. 
  At $E$ the diagonal and decomposition are no longer unique, but
  each can be defined by taking the limit along~$\mtilde\alpha$.
  The limiting diagonal $e$ in~$E$ then develops to a semicircle 
  contained in the equator, 
  bounded by the antipodal vertices $v_i$ and~$v_j$ of~$S$.

  We now define $\lambda$ so that it lifts to a path in $\Dk^*$ 
  which removes a slit sphere whose boundary contains the diagonal~$e$.
  Of the pieces adjacent to~$e$ at least
  one---call it~$P$---is isometric to a slit sphere,
  because there is only one piece isometric to a hemisphere 
  in the decomposition.
  Now consider the vertices developing into the image of~$e$.  
  Note that by the definition of~${\mathcal E}_k$, these vertices are either
  all the vertices with index between $i$ and~$j$ or all the complementary ones.
  But~$e$ decomposes the metric into two connected halves, and these are
  exactly the vertices that lie to one side of~$e$.
  In order for the lift of~$\lambda$ to remove~$P$,
  we revolve the vertices which develop into~$e$ 
  by an angle of~$2\pi$ about the fixed endpoints $v_i,v_j$;
  we must revolve in the direction of~$P$.
  By a slight perturbation of $\alpha$ and~$S$, if necessary, 
  we can assume that the orbits of the revolving vertices
  stay disjoint from the fixed vertices.  
  This ensures that $\lambda$ lies in~$\Tk^*$ and so has a unique lift 
  $\tilde\lambda$ in~$\Dk^*$.  
  
  The complement in~$E$ of the slit sphere~$P$ has two components.
  We can explicitly define the path~$\tilde\lambda$ by attaching these
  components to a path of lunes whose angles decrease from $2\pi$ to~$0$.
  By leaving the rest of the decomposition unchanged,
  we obtain a natural decomposition of the metrics along~$\tilde\lambda$.
  The decomposition of the terminal metric~$E'$ of~$\tilde\lambda$
  then agrees with that of~$E$, except that the one slit sphere
  is removed.  
  
  There is a unique lift $\mtilde\alpha^{-1}$ starting from $E'$;
  it terminates at some $D'\in F$.
  As with~$\mtilde\alpha$, the metrics along $\mtilde\alpha^{-1}$
  (aside from~$E'$) have a unique decomposition.
  This path of decompositions agrees with that along the reversed
  path~$\mtilde\alpha$ except one slit sphere has been removed;
  the limiting decomposition of this path therefore 
  agrees with the decomposition at~$E'$.
  The small truncation of~$D'$ has one less slit sphere than the 
  small truncation of~$D$, so $m(D')=m(D)-1$.
\end{proof}

\subsection{The space $\D_4$}

We now proceed to analyze in detail the structure of the
space~$\D_4$ of $4$-point metrics, proving that it is
diffeomorphic to a $5$-ball.
Recall that $\T_4^*$ is the $5$-manifold of all pairwise distinct 
spherical $4$-tuples $(p_1,p_2,p_3,p_4)$ modulo rotation;
the $5$-manifold $\T_4$ additionally includes the $4$-tuples
where either $p_1=p_3$ or $p_2=p_4$.

\begin{lemma}
  The space $\T_4^*$ is a $\C$-bundle
  whose base is the open $3$-ball with two proper arcs removed;
  it thus has the homotopy type of $\S^1\vee\S^1$, a figure-8.
  The space $\T_4$ is also a $\C$-bundle; the base is an open $3$-ball 
  with a proper arc and a point removed, so
  the homotopy type is $\S^1\vee\S^2$.
\end{lemma}

\begin{proof}
  Let us identify $\S^2$ with $\hat\C:=\C\cup\infty$.  
  We use a rotation to move $p_4$ to~$\infty\in\hat\C$.
  Then
  $$ 
    \T_4 \isom \Big(\{(p_1,p_2,p_3):
      p_1,p_3\in\C, p_2\in\hat\C, p_1\ne p_2\ne p_3\}
      \setm \{(p,\infty,p)\} \Big) \,\Big/\; \S^1, 
  $$
  where $\S^1$ acts diagonally by multiplication by~$e^{i\theta}$.

  We define $\X_4\subset\T_4$ to be the subset where $p_3=0$.
  Then
  $$
    \pi\colon \T_4\to\X_4,\qquad
    (p_1,p_2,p_3)\mapsto(p_1-p_3,p_2-p_3,0)
  $$
  is well-defined and indeed a deformation retract.
  In fact, $\pi$ is the desired bundle with fiber~$\C$:
  the $\C$-valued expression $p_3/(p_2-p_3)$ is $\S^1$-invariant,
  hence well-defined on~$\T_4$, and it gives a trivialization
  away from $p_2=\infty$;
  by symmetry a similar trivialization exists away from $p_1=p_3$.

  In our picture, $\T_4^*\subset\T_4$ is the subspace
  where $p_1\ne p_3$ and $p_2\ne\infty$.  So if we define
  $$
    \X_4^*:=
      \big\{(p_1,p_2,0): p_1, p_2\in\C\setm\{0\}, p_1\ne p_2\big\} / \S^1,
  $$
  the bundle $\T_4\to\X_4$ restricts to a $\C$-bundle $\T_4^*\to\X_4^*$.

  Now consider the smooth map $\phi\colon\X_4\to\C\times\R$ defined by
  $\phi(p_1,p_2,0) := (z,t)$, where
  $$
    z := \frac{p_1}{p_2}, \qquad
    t := \log\big(1+|p_1|\big)-\log\Big(1+\frac1{|p_2|}\Big)
      = \log\frac{1+|p_1|}{1+\sfrac{|z|}{|p_1|}}. 
  $$
  The map $\phi$ is well-defined since it is unchanged when multiplying
  $p_1$ and~$p_2$ by the same factor~$e^{i\theta}$.  
  Note that $z$ is a cross-ratio of the $4$-tuple
  $$(p_1,p_2,p_3,p_4)=(p_1,p_2,0,\infty).$$
  
  We now examine the behavior of $\phi$ on~$\X_4^*$; 
  here $z$ takes every value except~$0$ and~$1$.  
  For any fixed $z=\sfrac{p_1}{p_2}$, the function
  $$
    r\mapsto \log \frac{1+r}{1+\sfrac{|z|}r}
  $$
  is strictly increasing from $\R^+$ onto~$\R$.
  Thus $\X_4^*$ is diffeomorphic to its image
  $\phi(\X_4^*)=\big(\C\setm\{0,1\}\big)\times\R$,
  and so has the homotopy type of $\S^1\vee\S^1$.
  This completes the proof for~$\T_4^*$.

  The space $\X_4$ includes two open arcs in addition to~$\X_4^*$,
  namely the arcs~$\gamma_{13}$, where $p_1=p_3$,
  and~$\gamma_{24}$, where $p_2=p_4$.  Under~$\phi$, both
  arcs map to the line $z=0$; the coordinate $t$
  maps $(0,p_2,0)\in\gamma_{13}$ to $-\log\big(1+\tfrac1{|p_2|}\big)$, 
  and $(p_1,\infty,0)\in\gamma_{24}$ to $\log(1+|p_1|)$.
  Hence this coordinate
  maps $\R^+$ diffeomorphically onto~$\R^-$ and~$\R^+$, respectively.
  Thus $\phi(\X_4)$ is the union of $\phi(\X_4^*)$ with $0\times(\R\setminus 0)$.
  We conclude $\phi$ is a diffeomorphism from~$\X_4$ to its image,
  which is all of~$\C\times\R$ except for the origin 
  and the line $\{1\}\times\R$.
  This is homotopy-equivalent to $\S^2\vee\S^1$,
  completing the proof for~$\T_4$.
\end{proof}

Remembering that $\D_4^*$ is a covering space of $\T_4^*$,
we know it will be a $\C$-bundle over the corresponding
cover $\Y_4^*$ of $\X_4^*$.  Letting $\Zodd := 2\Z+1$ denote the
odd integers, we have:

\begin{lemma}
  The space $\D_4^*$ is a $\C$-bundle over $\Y_4^*\isom(\C\setm\Zodd)\times\R$.
\end{lemma}
\begin{proof}
  The fundamental group of $\T_4^*$ is $\Z *\Z$.
  The first generator~$\beta$ is a loop represented by $p_1$ circling
  counterclockwise around~$p_2$;
  in terms of $\phi(\X_4^*)=\big(\C\setm\{0,1\}\big)\times\R$,
  it corresponds to a loop around the deleted line at~$1$.
  The second generator~$\gamma$ is represented
  by $p_1$ circling clockwise around~$p_3$
  (or equivalently, counterclockwise around both $p_2$ and $p_4$);
  in $\phi(\X_4^*)$, it is a loop around the deleted line at~$0$.

  We claim that every lift of $\gamma$ to $\D_4^*$ is a loop,
  while the lifts of $\beta$ are not.  
  This means that the covering $\D_4^*\to\T_4^*$ is normal,
  corresponding to the normal subgroup generated by~$\gamma$.  
  The covering map $\Y_4^*\to\X_4^*$ can be then explicitly given as
  $(w,t)\mapsto(1+e^{\pi i w},t)$ with $\Y_4^*\isom(\C\setm\Zodd)\times\R$,
  and we can represent $\D_4^*$ as the pullback $\C$-bundle over~$\Y_4^*$.

  To prove the claim, we will work over the deformation retract~$\X_4^*$.
  Because $p_3=0$ and $p_4=\infty$ are antipodal, 
  truncations of metrics projecting to~$\X_4^*$ are not unique.
  To overcome this difficulty we adopt the convention that we always truncate 
  that end along an arc developing to the positive reals.  
  This corresponds to viewing the point $(p_1,p_2,0)$
  as a limit of $(p_1,p_2,\epsilon)$ for $\epsilon$ positive real.

  We fix a basepoint $P:=(1,2,0)\in\X_4^*\subset\T_4^*$.
  We can use our decomposition results to analyze any metric
  in the fiber over~$P$.  Because of the order of the four points in
  the circle, and our convention for drawing diagonals~$p_3p_4$,
  each decomposition piece is isometric to a slit sphere.
  Furthermore exactly two pieces are $3$-gons, with vertices
  either $p_1p_2p_3$ and $p_1p_3p_4$, or $p_1p_2p_4$ and $p_2p_3p_4$.
  Metrics in the fiber over~$P$ are completely characterized by
  the choice of one of these two types,
  together with the nonnegative number of slit spheres between the two
  triangles.

  More explicitly, suppose we truncate a $4$-point metric~$D$ in the
  fiber over~$P$, using again a diagonal $p_3p_4$
  which projects to the positive reals. 
  The $4$-gon $D'$ obtained has four interior angles $\alpha_i:=\alpha(p_i)$; 
  each of these is an integer multiple of~$2\pi$.
  The alternating sum
  $$
    \alpha_1-\alpha_2+\alpha_3-\alpha_4=2\pi n
  $$
  is unchanged if we add $2\pi$ to two consecutive angles,
  and hence is independent of the truncation $D'$ we have chosen.
  Note that the area of~$D'$ and the number $2(\alpha_2+\alpha_4)$ 
  are each multiples of~$4\pi$.
  Thus the Gauss--Bonnet Theorem shows that $n$ must be odd;
  the odd integer~$n$ uniquely characterizes a metric $D=D_n$ 
  in the fiber over~$P$.
  This also follows from our combinatorial picture:
  indeed, $D_n$ has a diagonal $p_1p_3$ when $n>0$ and 
  a diagonal $p_2p_4$ when $n<0$; the number of slit spheres
  along this diagonal is $\sfrac{(|n|-1)}2$.

  Now choose a truncation $D'$ of~$D_n$ with at
  least one slit sphere added at each end.
  The loops $\beta$ and~$\gamma$ in~$\X_4^*$
  can be lifted to paths in~$\D_4$ with initial point~$D_n$.
  These paths also have a unique lift starting at~$D'$ in
  the space $\P_4$ of $4$-gons, provided no angle~$\alpha_i$ becomes zero.  
  Indeed, this is easy to see if the paths are chosen to avoid antipodal 
  pairs in the $4$-tuple: then the unique geodesics between them
  lead to a unique continuation of the truncation.

  Let us now observe the effect of these lifts of $\beta$ and~$\gamma$
  on the angles~$\alpha_i$.
  Since only $p_1$ moves and the edge $p_3p_4$ stays 
  fixed (by our convention), the angle~$\alpha_3$ remains constant.
  For $\gamma$, as $p_1$ moves around $p_3$ (which is antipodal to~$p_4$)
  the edge $\overline{p_1p_4}$ rotates once around the sphere.
  This means that $\alpha_1$ and~$\alpha_4$ each decrease by~$2\pi$.
  Since the angle~$\alpha_2$ stays nearly constant all along~$\gamma$,
  it must return to its initial value.  We conclude that $n$ is fixed,
  and so $\gamma$ lifts to a loop in~$\D_4^*$.
  For~$\beta$, on the other hand, $\alpha_4$ is nearly constant,
  $\alpha_2$ increases by $2\pi$ and~$\alpha_1$ decreases by~$2\pi$;
  thus $n$ decreases by~$2$.  
  Thus $\beta$ changes the metric $D_n$ to~$D_{n-2}$.  
  This completes the proof of the claim.
\end{proof}

\begin{proposition}
  The space~$\D_4$ is diffeomorphic to an open 5-ball.
\end{proposition}
\begin{proof}
  Remember that $\X_4\setm\X_4^*$ consists of two arcs $\gamma_{13}$
  and~$\gamma_{24}$.  The decision of where to fill in copies of these
  arcs into the covering space $\Y_4^*$ is purely combinatorial, based
  on the decompositions of the metrics in $\D_4\setm\D_4^*$.

  We can connect our basepoint $P=(1,2,0)$
  to the fibers over these arcs by simple paths:
  the first path moves only~$p_1$ (along the reals to~$0$); 
  the second moves only~$p_2$ (along the reals to~$\infty$).  
  We use these paths to extend the odd labels~$n$ from the fiber over~$P$.   
  Now we note that $p_i=p_j$ is allowed in $\D_4$ exactly when there
  is no diagonal $p_ip_j$.
  Thus the arc $\gamma_{13}$ gets filled in exactly at the levels
  with negative label $n<0$.  
  Similarly the arc $\gamma_{24}$ gets filled in exactly 
  at the levels with $n>0$.

  This means that at each level, exactly one of the two arcs gets filled in.
  Again, the decision is purely combinatorial: whenever we fill in a single
  point, we fill in an entire arc into~$\Y_4^*$,
  and $\D_4$ contains the entire $\C$-bundle over this arc.
  In the picture $\Y_4^*\isom(\C\setm\Zodd)\times\R$,
  we fill in the lower half ($t<0$) of the missing lines where $z=n<0$,
  and the upper half ($t>0$) of the others where $z=n>0$; 
  this leaves us with a space $\Y_4$ diffeomorphic to a three-ball.  
  Finally, $\D_4$ is a $\C$-bundle over this ball, 
  which must be trivial, meaning that $\D_4$ is a five-ball.
\end{proof}

\section{Injectivity and Properness}

To show the classifying map $\Phi$ is continuous, injective, and proper,
we follow the approach we used for the case $k=3$ of triunduloids.
The proof of properness is where substantial modifications are needed.
In particular, since the classifying space is~$\Dk$, not~$\Tk$,
we need to characterize the compact subsets of~$\Dk$ for general~$k$.

\subsection{Injectivity}

As for triunduloids we have for coplanar $k$-unduloids of genus zero:
\begin{theorem}\label{th:inj}
  The classifying map $\Phi$ is injective.
\end{theorem}
\begin{proof}
  The proof is essentially the same as for triunduloids.
  We sketch the main steps and refer to \cite{gks} for details.
 
  Suppose two $k$-unduloids $M,M'\in\Mk$ map to the same 
  $k$-point metric $D:=\Phi(M)=\Phi(M')\in\Dk$.  We need to show $M=M'$. 
  We pull back the $\k$-Hopf bundle $\S^3\to\S^2$ to 
  the disk $D$ by the developing map.  Taking the universal cover of
  the fibers, we obtain a line bundle $E\to D$.
  It carries a Riemannian structure, locally isometric to~$\S^3$. 
  The $\k$-Hopf flow acts isometrically on~$E$, moving points vertically
  along the fibers.

  The cousin disks $\plustilde M$ and $\plustilde{M'}$ can be considered 
  sections of~$E$.   Over each of the $k$ rays of the decomposition of~$D$
  given by Theorem~\ref{thm:decomp}, the $\k$-Hopf flow foliates the 
  bundle~$E$ with unduloid cousins.  
  Over each ray, the cousin disks are asymptotic to one of the leaves.
  But across each end of~$M$ and~$M'$ the \emph{periods}
  (see \cite[Sect.~3.3]{gks}) vanish,
  which by \cite[Lem.~3.7]{gks} implies for the cousins that 
  the $\k$-Hopf flow can move $\plustilde M$ to become
  asymptotic to~$\plustilde{M'}$ over all rays simultanuously.  
  This allows us to deal with
  $\plustilde M$ and $\plustilde {M'}$ as if they were compact;
  as in \cite[Thm.~3.8]{gks} the 
  maximum principle then shows 
  $\plustilde M=\plustilde {M'}$.
  Therefore $M=M'$, up to translation.
\end{proof}

\subsection{Properness}

To show properness, we want to characterize compact subsets of~$\Dk$.
We need the following notion.
The \emph{vertex distance} $d(D)$ of $D\in\Dk$ is the minimum  
distance between any pair of vertices in the completion $\hat D$.

\begin{proposition} \label{pr:cptK}
  For any compact $\K\subset\Dk$, there exist $m,\epsilon>0$ 
  such that for all $D\in\K$\\
  \ia\ the vertex distance $d(D)$ is at least~$\epsilon$, and \\
  \ii\ there is a small truncation $\Dbox$ with at 
     most~$m$ decomposition pieces. 
\end{proposition}
\begin{proof}
  \ia\  
  The distance between any two given vertices
  depends continuously on $D\in\Dk$ and is always positive;
  thus it assumes a positive minimum over the compact set~$K$.
  Minimizing over the finitely many pairs of distinct vertices
  still yields a positive minimum for~$d(D)$.
  (It is not hard to see that, in fact, $d(D)$ is realized
  by the spherical distance of some pair of distinct
  points of~$\Pi(D)\in\Tk$.)

  \smallskip\noindent\ii\
  Consider a metric $D\in\Dk$ which admits a small truncation
  with $n$ decomposition pieces.  
  If no piece is a lune or a concave polygon with angle~$\pi$ then 
  $D$ has a neighborhood of metrics with the same number of pieces.
  Otherwise, at most $k$~additional pieces are needed.

  Now let $V_j\subset\Dk$ be the set of metrics
  which admit a small truncation with at most $j-k$ decomposition pieces.  
  The above discussion shows $V_j$ has an 
  open neighborhood $U_j\subset\Dk$, such that any metric in~$U_j$ 
  admits a truncation with no more than $j$~decomposition pieces. 
  The $\{U_j\}$ define an open exhaustion of~$\Dk$, so that the
  compact set $K$ is contained in some~$U_m$.
%
\end{proof}

\begin{lemma}\label{le:curv}
  Suppose $\K$ is a subset of $\Dk$ such that the vertex
  distance of each $D\in\K$ is at least~$\epsilon>0$.
  Then there is $C=C(\epsilon)>0$ such that the curvature estimate
  $\sup_M |A|<C(\epsilon)$ holds for each $M\in\Phi^{-1}(\K)$.
\end{lemma}
\begin{proof}
  We modify the argument given in \cite[Sect.~4]{gks}.
  Suppose on the contrary that $\sup_{M_j}|A|\to\infty$ with $\Phi(M_j)\in \K$.
  Then by \cite[Lem.~4.2]{gks} we find a sequence of rescaled surfaces~$M'_j$,
  converging to some minimal surface~$N$ which is Alexandrov-symmetric
  across the plane~$P$.
  Note that, as in \cite[Prop.~4.4]{gks}, the $M_j$ have genus zero,
  meaning that $M_j\cap P$ has no closed loops.  Since the convergence
  is with multiplicity one, each component of $N\cap P$ is likewise
  an unbounded curve.

  By \cite[Lem.~4.3]{gks} the minimal surface~$N$ contains a simple
  closed curve~$\gamma$ with nonzero force, symmetric with respect to~$P$
  and not contained in~$P$.  So the surfaces~$M'_j$ contain closed 
  curves~$\gamma'_j$, whose lengths and forces converge to those of~$\gamma$.  
  Now consider the closed curves~$\gamma_j$ on the unrescaled surfaces~$M_j$.
  Note that the length of~$\gamma_j$ converges to~$0$.
  For sufficiently large~$j$, however, the force of $\gamma'_j$ 
  is nonzero, and so $\gamma'_j$ as well as~$\gamma_j$ 
  must be homologically nontrivial.  
  But since the~$M_j$ have genus zero this means that 
  the curve~$\gamma_j$ must connect two different components of $M_j\cap P$.
  Therefore, the conjugate curve $\tilde\gamma^+_j$
  to its \uhalf/ $\gamma^+_j$ Hopf projects to a curve~$\eta_j$ 
  in $\Phi(M_j)\in\Dk$ which connects two different vertices.  

  By definition, the length of any such curve is at least the vertex 
  distance $d(D)\ge\epsilon$.
  Since the Hopf projection $\Pik$ has Lipschitz constant~$2$,
  it follows that the length of $\tilde\gamma^+_j$, and thus also
  the length of $\gamma_j$, is bounded below.
  This contradicts the fact that $\gamma_j$ has length converging to~$0$.
\end{proof}

\begin{lemma}\label{le:ball}
  Given a compact subset $\K\subset\Dk$ there is $R>0$ with the following
  property.  For each $M\in \Phi^{-1}(\K)$ there exists a ball $B_R\in\R^3$
  of radius~$R$ such that
  the noncompact components of $M\setm B_R$ are exactly $k$ annular ends, 
  each of which is contained in a solid cylinder of radius~$3$.
\end{lemma}
\begin{proof}
  The enclosure theorem, proved in~\cite{kk} for embedded surfaces,
  carries over without change to say that any \alex-emb/ \CMC/
  surface $M$ is contained
  in the union of $k$~solid cylindrical rays as well as a finite
  number of balls of radius $11$ and solid cylindrical segments of radius~$3$; 
  the numbers of these depend only on~$k$ and the genus.
  Collecting the balls and finite cylindrical segments into one larger ball,
  we see that for each $M\in\Mk$ there is an $R=R(M)$ with the property
  desired.  Note that $\Mk$ contains sequences of $k$-unduloids 
  where the length of some midsegment, and hence~$R$, diverges.
  However we claim that, over~$\K$, there is a uniform choice of~$R$.

  Suppose, by contradiction, that for some
  sequence $M_j\in \Phi^{-1}(\K)$ we have $R(M_j)\to\infty$.
  Consequently, the enclosing sets for the~$M_j$ 
  have solid cylindrical segments with
  length~$l(j)$ tending to infinity.
  These segments contain annuli~$A_j\subset M_j$.
  By the uniform bound on the vertex distance from 
  Proposition~\ref{pr:cptK}\ia,
  the force carried by these annuli are uniformly bounded away from $0$.
  Thus, as $l(j)\to\infty$, it follows as in~\cite[Thm.~5.1]{kks} 
  that after 
  rigid motion a subsequence of the~$A_j$ 
  converge to some (complete) unduloid.

  Since this unduloid has a cousin which Hopf-projects to
  a line of slit spheres, we find that the spherical metrics
  for~$A_j$ must contain increasing numbers of slit spheres.
  This contradicts Proposition~\ref{pr:cptK}\ii.
\end{proof}

\begin{theorem}                      \label{th:proper}
  The classifying map $\Phi\colon\Mk\to\Dk$ is continuous and proper.
\end{theorem}
\begin{proof}
  Given a compact $\K\subset\Dk$, Lemma~\ref{le:curv} (combined with
  Proposition~\ref{pr:cptK}\ia) gives a uniform
  curvature estimate for $k$-unduloids in $\Phi^{-1}(\K)$.  The continuity 
  of~$\Phi$ on this set then follows just as in~\cite[Thm.~4.5]{gks}.

  To show that $\Phi^{-1}(\K)$ is compact,
  we now prove that any sequence $M_j\in \Phi^{-1}(\K)$
  subconverges with respect to the topology on~$\Mk$.
  The~$M_j$ are given only up to horizontal translation; to
  get convergence, we choose representatives such that
  the enclosing balls $B_R$ from Lemma~\ref{le:ball}
  are centered at the origin.
  Note also that by~\cite{kk}, there is a uniform area bound, 
  for any coplanar genus-zero $k$-unduloid $M$,
  of the form $\area(M_j\cap B_r)\le Cr^2$.
  This area bound, together with the curvature estimate
  of Lemma~\ref{le:curv}, now gives a subsequence of the~$M_j$
  that converges (on each compact subset) to some \CMC/ surface~$M$.

  By smooth convergence, it is clear that $M$ is \alex-emb/.
  Because a sequence of contractible loops has a contractible limit,
  $M$ has genus zero and at most $k$ ends.
  On the other hand, by Lemma~\ref{le:ball}, on each $M_j$
  there are $k$ closed curves in $B_{R+1}(0) \setm B_{R}(0)$
  which bound the $k$ ends.  These $k$ sequences of
  curves subconverge to $k$ closed curves which are disjoint in
  the limit surface $M$, and thus bound $k$ different ends.
  Thus $M$ is a coplanar $k$-unduloid, as desired.
\end{proof}

\section{Surjectivity of the classifying map}

\subsection{Moduli spaces of coplanar $k$-unduloids}

Given a minimal or \CMC/ surface~$M$, a solution~$u$ to the
linearized mean-curvature equation
$L(u):=\Delta_M u+|A|^2 u=0$ is called a \emph{Jacobi field}.
The surface~$M$ is called 
\emph{nondegenerate} if all square-integrable Jacobi fields
$u\in L^2(M)$ vanish.

Let $\M_{g,k}$ denote the moduli space of all $k$-unduloids of genus~$g$.
Kusner, Mazzeo and Pollack~\cite{kmp} showed that~$\M_{g,k}$
is a real-analytic variety (for $k\ge3$), and that near any
nondegenerate surface it is a $(3k-6)$--dimensional manifold.
We now prove a similar result for coplanar $k$-unduloids:

\begin{theorem} \label{th:kmp}
  The moduli space $\Mk$ of coplanar $k$-unduloids of genus zero
  is locally a real-analytic variety of finite dimension.
  In a neighborhood of a nondegenerate coplanar $k$-unduloid,
  it is a manifold of dimension~$2k-3$.
\end{theorem}
In their proof, \cite{kmp} consider the \emph{premoduli space} 
$\wtilde\M_{g,k}$, which is the space of $k$-unduloids of genus~$g$
\emph{before} dividing by euclidean motions.
At a nondegenerate surface~$M$, they show the premoduli space $\wtilde\M_{g,k}$
is locally a $3k$-dimensional Lagrangian submanifold of a symplectic vector 
space~$W$.  Here, each end of~$M$ has a
$6$-dimensional space of geometric Jacobi fields, 
induced by translations, rotations, and change of asymptotic necksize.
The direct sum of these, over the $k$~ends of~$M$, forms the 
$6k$-dimensional space~$W$.
The tangent space to $\wtilde\M_{g,k}$ at~$M$ is the Lagrangian subspace 
consisting of all globally defined Jacobi fields on~$M$.  
\begin{proof}
  In the premoduli space $\wtilde\M_{0,k}$ of all genus-zero $k$-unduloids,
  let $\wtilde\M_k$ denote the subset consisting 
  of coplanar $k$-unduloids with Alexandrov symmetry
  across the $xy$-plane~$P$.
  The reflection $\tau$ across~$P$ acts on~$\wtilde\M_{0,k}$,
  fixing~$\wtilde\M_k$.

  Near any surface~$M$, the variety~$\wtilde\M_{0,k}$
  is defined by real-analytic 
  equations on a finite-dimensional vector space $W'\supset W$.
  Our premoduli space~$\wtilde\M_k$ is obtained locally by adding one more set
  of linear equations, $u\after\tau=u$, and is thus still real-analytic.
  Dividing by euclidean motions of the $xy$-plane~$P$, we find the same
  is true for~$\Mk$.
 
  To compute the dimension of~$\Mk$ we show that 
  $\wtilde\M_k$ has dimension~$2k$ in the neighborhood of a
  nondegenerate surface~$M$.  
  We use the following transversality argument. 
  The reflection~$\tau$ induces an automorphism $\tau^*$
  on the $6k$-dimensional vector space~$W$ which encodes the asymptotic
  behavior of Jacobi fields on~$M$. 
  This automorphism preserves the natural symplectic form on~$W$.  
  The global Jacobi fields on~$M$ form a Lagrangian subspace $L\subset W$,
  which corresponds to the tangent space to~$\wtilde\M_{0,k}$ at~$M$.
  Since $\tau$ takes one global Jacobi field on~$M$ to another,
  $\tau^*$ carries $L$ into itself.  
  The fixed point set of~$\tau^*$ is a $4k$-dimensional vector space~$V$, 
  consisting of those asymptotic changes that keep the ends in the plane~$P$.
  The following lemma then says that $L\cap V$, 
  corresponding to the tangent space of $\wtilde\M_k$ at~$M$,
  is $2k$-dimensional.
\end{proof}

\begin{lemma}
  Suppose $L$ is a Lagrangian (maximal isotropic) subspace of
  a symplectic vector space $W$, and suppose $G$ is a compact group
  of symplectic automorphisms of~$W$ preserving~$L$.
  Then the fixed-point set $W^G$ of~$G$ is a symplectic subspace 
  $W^G\subset W$, and $L\cap W^G$ is Lagrangian in $W^G$.
\end{lemma}
\begin{proof}
  By averaging over the group~$G$, we can choose a compatible almost complex
  structure~$J$ on~$W$ that is also invariant under~$G$.
  That is, $g\after J=J\after g$ for all $g\in G$.
  Let $W^g := \ker(g-I)$, which is preserved by $J$.
  Then the fixed set $W^G$ is the intersection $W^G=\bigcap_{g\in G} W^g$.
  In particular, $J(W^G)=W^G$, showing this is symplectic.

  Because $L$ is Lagrangian, we have $W=L\oplus JL$, so
  any vector in $W$ can be represented as $v+Jw$, where $v,w\in L$.
  If, moreover, $v+Jw$ is in $W^g$ 
  then $g(v+Jw)=v+Jw$ and so $-(g-I)v = (g-I)Jw = J(g-I)w$.
  But since $g$ (and therefore $g-I$) preserves~$L$, 
  this equates a vector in~$L$ with one in~$JL$.  
  Thus $(g-I)v=0=(g-I)w$, so that in fact $v,w\in L\cap W^g$.
  Since the preceding argument holds for all $g\in G$ we obtain 
  a decomposition $W^G=(L\cap W^G)\oplus J(L\cap W^G)$, 
  meaning that $L\cap W^G$ is Lagrangian in~$W^G$.  
\end{proof}

\subsection{Existence of nondegenerate coplanar $k$-unduloids}

A \emph{minimal $k$-noid} is a finite-total-curvature minimal surface
of finite topology with genus zero which has only catenoidal ends.
Jorge and Meeks determined the $k$-noids which are coplanar and dihedrally 
symmetric.  They have only two umbilic points, so that the result of
Montiel and Ros~\cite{mr} applies (in the form mentioned in~\cite{gks}) to give:
\begin{proposition} \label{pr:knoids}
  The Jorge--Meeks $k$-noids are nondegenerate.
\qed
\end{proposition}

As in the case of triunduloids, we want to use these minimal surfaces
in the gluing technique of Mazzeo and Pacard~\cite{mp} to produce nondegenerate
coplanar \CMC/ surfaces.  Their theorem generalizes to an
equivariant setting as follows:

\begin{theorem}\label{th:mp}
  Let $M_0\subset\R^3$ be a nondegenerate minimal $k$-noid that
  can be oriented so that each end has normal pointing inwards towards its
  asymptotic axis.  Assume $M_0$ has symmetry group~$G$ in~$\R^3$,
  consisting of isometries preserving the normal.
  Then for some $\epsilon_0$,
  there is a family $M_{\epsilon}$, $0<\epsilon<\epsilon_0$,
  of nondegenerate \CMC/ surfaces with $k$~embedded ends,
  also with symmetry group~$G$, such that on any compact set in~$\R^3$,
  the dilated surfaces
  $\tfrac 1{\epsilon} M_{\epsilon}$ converge smoothly and uniformly to~$M_0$.
\end{theorem}

\begin{proof}[Sketch of proof]
The proof of this theorem follows exactly the argument
of Mazzeo and Pacard~\cite{mp}, so we merely outline
how to adapt their argument to our $G$-invariant setting.
They first glue exactly unduloidal ends to~$M_0$, 
truncated along $k$~circles, and then consider normal perturbations.
They slightly perturb the bounding circles, and
derive solution operators for the respective Dirchlet problems
for the mean curvature equation, both over~$M_0$ and over the ends.  
Unique solutions to these Dirichlet problems are established
as fixed points of a contraction mapping \cite[Cor.~5,~8]{mp};
thus $G$-invariant Dirichlet data lead to $G$-invariant solutions. 
The surface $M_\epsilon$ is then obtained by a Cauchy matching
of a pair of Dirichlet solutions: 
This matching is accomplished \cite[Prop.~32]{mp}
by applying the Leray-Schauder fixed-point theorem
to a nonlinear operator~$F_\epsilon$ on the Dirichlet data.
The operator~$F_\epsilon$ is designed to have a fixed point exactly
when the first-order Cauchy data match.

Given a $G$-invariant initial surface~$M_0$ 
we note that $G$ acts by composition on all relevant function spaces.
Also, the operators used in the various fixed-point arguments, being
geometrically defined, all commute with~$G$.
Therefore we can simply replace the function spaces by their
closed linear subspaces consisting of $G$-invariant functions.
Clearly, all estimates in~\cite{mp} continue to hold, 
and so the contraction-mapping and Leray--Schauder fixed-point theorems
still apply, giving us a $G$-invariant \CMC/ surface~$M_\epsilon$.

The arguments of \cite[Sect.~12]{mp} show that
for sufficiently small~$\epsilon$ any function in~$L^2(M)$
(whether $G$-invariant or not)
satisfying the linearized mean curvature equation vanishes.
Thus~$M_\epsilon$ is nondegenerate.
\end{proof}

Now let~$M_0$ be the Jorge--Meeks $k$-noid of Proposition~\ref{pr:knoids},
and $G$ be its symmetry group.
Applying Theorem~\ref{th:mp} gives:
\begin{corollary}\label{co:nondegkund}
  For each $k\ge 3$ there exists a nondegenerate
  coplanar genus-zero $k$-unduloid.
\qed
\end{corollary}

\subsection{Proof of Main Theorem and Corollary}

In \cite[Sect.~5.3]{gks} we show:
\begin{theorem} \label{th:homeo}
  If $f\takes X\to Y$ is a continuous, proper, injective map from
  a $d$-dimensional real-analytic variety $X$ to a connected $d$-manifold $Y$,
  then $f$ is surjective and thus a homeomorphism. \qed
\end{theorem}

We now have all the ingredients needed to prove:
\begin{theorem}
  The classifying map $\Phi\colon \Mk\to\Dk$ is a homeomorphism
  between connected $(2k-3)$--manifolds.
\end{theorem}
\begin{proof}
  To apply Theorem~\ref{th:homeo} to $\Phi$ we need to verify its assumptions.
  The space $\Dk$ is a connected $(2k-3)$--manifold by Theorems~\ref{th:mfld}
  and~\ref{th:connected}.  The map $\Phi$ is proper and injective by
  Theorems~\ref{th:proper} and~\ref{th:inj}.

  From Theorem~\ref{th:kmp}, the moduli space~$\Mk$ of coplanar $k$-unduloids
  is locally a real-analytic variety; injectivity of $\Phi$ shows it has
  dimension at most $2k-3$.
  But by Corollary~\ref{co:nondegkund} there is a nondegenerate coplanar
  $k$-unduloid, so then by Theorem~\ref{th:kmp}, $\Mk$ is
  a $(2k-3)$--manifold in some neighborhood of this $k$-unduloid.

  Thus $\Mk$ has dimension exactly $2k-3$, so
  Theorem~\ref{th:homeo} applies to conclude that $\Phi$ is a homeomorphism.
\end{proof}

Together with Theorem~\ref{th:necksizes}, this proves our Main Theorem.

The corollary on necksizes stated in the introduction collects the
statements of Corollary~\ref{co:oddk} and Corollary~\ref{co:evenk}.
The existence statement follows 
from the example of Section~\ref{se:polmet}; in particular
the example shows that \eqref{necksizesumodd} and \eqref{necksizesumeven}
are sharp.


\end{document}